\documentclass[12pt,twoside]{article}

\begin{document}

\centerline{}

\centerline {\Large{\bf On soft linear spaces and soft normed linear spaces}}

\centerline{}

\newcommand{\mvec}[1]{\mbox{\bfseries\itshape #1}}

\centerline{\bf {Sujoy Das$^1$, Pinaki Majumdar$^2$ and  S.K. Samanta$^3$}}
\centerline{}

\centerline{$^1$Department of Mathematics,}
\centerline{Bidhan Chandra College}\centerline{Asansol-713304, West Bengal, India }
\centerline{e-mail: sujoy\_math@yahoo.co.in}
\centerline{}
\centerline{$^2$Department of Mathematics,}
\centerline{MUC Women's College}\centerline{Burdwan-713104, West Bengal, India }
\centerline{e-mail: pmajumdar2@rediffmail.com}
\centerline{}
\centerline{$^3$Department of Mathematics}\centerline{ Visva Bharati}\centerline{Santiniketan, West Bengal, India.}
\centerline{e-mail: syamal\_123@yahoo.co.in}
\centerline{}

\newtheorem{Theorem}{\quad Theorem}[section]

\newtheorem{definition}[Theorem]{\quad Definition}

\newtheorem{theorem}[Theorem]{\quad Theorem}

\newtheorem{remark}[Theorem]{\quad Remark}

\newtheorem{corollary}[Theorem]{\quad Corollary}

\newtheorem{note}[Theorem]{\quad Note}

\newtheorem{lemma}[Theorem]{\quad Lemma}

\newtheorem{example}[Theorem]{\quad Example}

\newtheorem{proposition}[Theorem]{\quad Proposition}

\begin{abstract}
{\emph{In this paper an idea of soft linear spaces and soft norm on soft linear spaces are given and some of their properties are studied. Soft vectors in soft linear spaces are introduced and their properties are studied. Completeness of soft normed linear space is also studied and equivalent soft norms and convex soft sets are studied in soft normed linear space settings.}}

\end{abstract}

{\bf Keywords:}  \emph{Soft sets, soft elements, soft vectors, soft linear spaces, soft normed linear spaces, soft Banach spaces, equivalent soft norm, convex soft sets.}\\

{\bf 2000 MSC No:} 08A02.

\section{Introduction}
In the year 1999, Molodtsov \cite{Mol} initiated the theory of soft sets as a new mathematical tool for dealing with uncertainties. He has shown several applications of this theory in solving many practical problems in economics, engineering, social science, medical science, etc. Research works in soft set theory and its applications in various fields have been progressing rapidly since Maji et al. (\cite{Maji1},\cite{Maji2}) introduced several operations on soft sets and applied it to decision making problems. In the line of reduction and addition of parameters of soft sets some works have been done by Chen \cite{Chen}, Pei and Miao \cite{PM} , Kong et al. \cite{Kong} , Zou and Xiao \cite{ZX}. Aktas and Cagman \cite{AC} introduced the notion of soft group and discussed various properties. Jun (\cite{Jun1},\cite{Jun2}) investigated soft BCK/BCI -- algebras and its application in ideal theory. Feng et al. \cite{F1} worked on soft semirings, soft ideals and idealistic soft semirings.  Ali et al. \cite{Ali} and Shabir and Irfan Ali (\cite{Ali},\cite{SI}) studied soft semigroups and soft ideals over a semi group which characterize generalized fuzzy ideals and fuzzy ideals with thresholds of a semigroup. The idea of soft topological spaces was first given by M. Shabir, M. Naz \cite{MM} and mappings between soft sets were described by P. Majumdar, S. K. Samanta \cite{Pin1}. Feng et al. \cite{F2} worked on soft sets combined with fuzzy sets and rough sets. A. Recently in (\cite{S1},\cite{S2}) we have introduced a notion of soft real sets, soft real numbers, soft complex sets, soft complex numbers and some of their basic properties have been investigated. Some applications of soft real sets and soft real numbers have been presented in real life problems. Two different notions of 'soft metric' are presented in (\cite{S3}, \cite{S14}) and some properties of soft metric spaces are studied in both cases.

\noindent In this paper we have introduced a notion of soft linear space and soft normed linear space. In section 2, some preliminary results are given. In section 3, a notion of `soft linear space' is given and various properties of soft linear spaces are studied. In section 4, a definition of `soft vector'in a soft linear space is given and various properties of soft vectors are studied in details with examples and counter examples. A notion of `soft norm'in a soft linear space is introduced in section 5. It has been shown that every `soft normed linear space' is also a `soft metric space'\cite{S3}. In that section, completeness of soft normed linear spaces, equivalent soft norms and convex soft sets are studied in soft normed linear space settings.

\section{Preliminaries}

\begin{definition}
\cite{Mol} Let $U$ be an universe and $E$ be a set of parameters. Let ${\mathcal P}(U)$ denote the power set of $U$ and $A$ be a non-empty subset of $E$. A pair $(F,A)$ is called a soft set over $U$, where $F$ is a mapping given by $F:A\to {\mathcal P}(U)$.  In other words, a soft set over $U$ is a parametrized family of subsets of the universe $U$. For $\varepsilon \in A,\ F(\varepsilon )$ may be considered as the set of $\varepsilon $ -- approximate elements of the soft set $(F,A)$.  
\end{definition}

\begin{definition}
\cite{F2} For two soft sets $(F,A)$ and $(G,B)$ over a common universe $U$, we say that $(F,A)$ is a soft subset of $(G,B)$ if 

\begin{enumerate}
\item  $A\subseteq B$ and 

\item  for all $e\in A,\ F(e)\subseteq $ $G(e)$.    We write $(F,A)\widetilde{\subset }$ $(G,B)$.
\end{enumerate}

\noindent $(F,A)$ is said to be a soft superset of $(G,B)$, if $(G,B)$ is a soft subset of $(F,A)$. We denote it by$(F,A)\widetilde{\supset }$ $(G,B)$.
\end{definition}

\begin{definition}
\cite{F2} Two soft sets $(F,A)$ and $(G,B)$ over a common universe $U$ are said to be equal if $(F,A)$ is a soft subset of $(G,B)$ and $(G,B)$ is a soft subset of $(F,A)$.
\end{definition}

\begin{definition}
\cite{F2} The complement of a soft set $(F,A)$ is denoted by 

\noindent ${(F,A)}^c=(F^c,A)$, where $F^c:A\to {\mathcal P}(U)$ is a mapping given by $F^c\left(\alpha \right)=U-F(\alpha )$, for all $\alpha \in A$.
\end{definition}

\begin{definition}
\cite{Maji2}  A soft set $(F,E)$ over $U$ is said to be an \textit{absolute} soft set denoted by${\rm \ }\check{U}$ if for all $\varepsilon \in E,\ \ F\left(\varepsilon \right)=U$. 
\end{definition}

\begin{definition}
\cite{Maji2}  A soft set $(F,E)$ over $U$ is said to be a \textit{null} soft set denoted by${\rm \ }\Phi $ if for all $\varepsilon \in E,\ F\left(\varepsilon \right)=\emptyset $. 
\end{definition}

\begin{definition}
\cite{Maji2} The union of two soft sets $(F,A)$ and $(G,B)$ over the common universe $U$ is the soft set 

\noindent $\left(H,C\right),$ where $C=A\cup B$ and for all $e\in C,$
\[H\left(e\right)=\left\{ \begin{array}{c}
F\left(e\right)\ \ \ \ \ \ \ \ \ \ \ \ \ \ \ \ \ \ \ \ \ \ \ \ if\ e\in A-B \\ 
G\left(e\right)\ \ \ \ \ \ \ \ \ \ \ \ \ \ \ \ \ \ \ \ \ \ \ \ \ if\ e\in B-A \\ 
F\left(e\right)\cup G\left(e\right)\ \ \ \ \ \ \ \ \ \ \ if\ e\in A\cap B. \end{array}
\right.\] 
We express it as $\left(F,A\right)\widetilde{\cup }\left(G,B\right)=(H,C)$.
\end{definition}

\noindent The following definition of intersection of two soft sets is given as that of the bi- intersection in \cite{F1}.

\begin{definition}
\cite{F1}  The intersection of two soft sets $(F,A)$ and $(G,B)$ over the common universe $U$ is the soft set $\left(H,C\right),$ where $C=A\cap B$ and for all $e\in C,$ $H\left(e\right)=F(e)\cap G(e)$. We write $\left(F,A\right)\widetilde{\cap }\left(G,B\right)=(H,C)$.
\end{definition}

\noindent Let $X$ be an initial universal set and $E$ be the non-empty set of parameters.

\begin{definition}
\cite{MM}  The difference $(H,E)$ of two soft sets $(F,E)$ and $(G,E)$ over $X$, denoted by $\left(F,E\right)\ \backslash (G,E)$, is defined by $H\left(e\right)=F(e)\backslash G(e)$ for all $e\in E$.
\end{definition}

\begin{proposition}
\label{p1}
\cite{MM} Let $(F,E)$ and $(G,E)$ be two soft sets over $X$. Then

\begin{itemize}  
\item [{(i).}] ${((F,E)\widetilde{\cup }(G,E))}^c={\left(F,E\right)}^c\widetilde{\cap }{\left(G,E\right)}^c$ 
\item [{(ii).}] ${((F,E)\widetilde{\cap }(G,E))}^c={\left(F,E\right)}^c\widetilde{\cup }{\left(G,E\right)}^c.$
\end{itemize} 
\end{proposition}

\begin{definition}
\cite{S1} Let $\ X\ be$ a non-empty set and $E\ $be a non-empty parameter set. Then a function

\noindent $\varepsilon :E\to X$  is said to be a soft element of $X.$ A soft element $\varepsilon $ of $X$ is said to belongs to a soft set $A$ of $X$, which is denoted by $\varepsilon \widetilde{\in }A,\ $if $\varepsilon \left(e\right)\in A\left(e\right),\ \forall e\in E.$ Thus for a soft set A of $X$ with respect to the index set $E,$ we have  $A\left(e\right)=\left\{\varepsilon \left(e\right),\varepsilon \widetilde{\in }A\right\},$ $e\in E.$
\end{definition}

\noindent It is to be noted that every singleton soft set (a soft set $(F,E)$ for which $F(e)$ is a singleton set, $\forall e\in E$) can be identified with a soft element by simply identifying the singleton set with the element that it contains $\forall e\in E.$

\begin{definition}
\cite{S1} Let $R$ be the set of real numbers and $\mathcal{B}(R)$ the collection of all non-empty bounded subsets of $R$ and $A$ taken as a set of parameters. Then a mapping  $F:A\to \mathcal{B}(R)$ is called a \textit{soft real set}. It is denoted by $(F,A)$. If specifically $(F,A)$ is a singleton soft set$,$ then after identifying $(F,A)$ with the corresponding soft element, it will be called a \textit{soft real number}.
\end{definition}

\noindent We use notations $\tilde{r},\ \tilde{s},\ \tilde{t}$ to denote soft real numbers whereas $\overline{r},\ \overline{s},\ \overline{t}$ will denote a particular type of soft real numbers such that $\overline{r}\left(\lambda \right)=r,\ $for all $\lambda \in A$ etc. For example $\overline{0}$ is the soft real number where $\overline{0}\left(\lambda \right)=0,\ $for all $\lambda \in A.$

\begin{definition} \cite{S3}
A mapping $d:SE(\check{X})\times SE(\check{X})\to {{\mathcal R}\left(A\right)}^*$ , is said to be a \textit{soft metric} on the soft set $\check{X}$ if $d$ satisfies the following conditions:

\begin{itemize}
\item [{(M1).}] $d\left(\tilde{x},\tilde{y}\right)\widetilde{\ge }\overline{0},\ $ for all $\tilde{x},\tilde{y}\widetilde{\in }\check{X}$.

\item [{(M2).}] $d\left(\tilde{x},\tilde{y}\right)=\ \overline{0}$ if and only if $\tilde{x}=\tilde{y}.$

\item [{(M3).}]  $d\left(\tilde{x},\tilde{y}\right)=d\left(\tilde{y},\tilde{x}\right)$ for all $\tilde{x},\tilde{y}\widetilde{\in }\check{X}$.

\item [{(M4).}] For all $\tilde{x},\tilde{y},\tilde{z}\widetilde{\in }\check{X},d\left(\tilde{x},\tilde{z}\right)\widetilde{\le }d\left(\tilde{x},\tilde{y}\right)+d\left(\tilde{y},\tilde{z}\right)$
\end{itemize}

\noindent The soft set $\check{X}$ with a soft metric $d$ on $\check{X}$ is said to be a \textit{soft metric space} and is denoted by $(\check{X},d,A)$ or $(\check{X},d)$. (M1), (M2), (M3) and (M4) are said to be soft metric axioms.
\end{definition} 

\begin{theorem}{(Decomposition Theorem)}  \cite{S3}
If a soft metric $d$ satisfies the condition:

\noindent (M5).For $\left(\xi ,\eta \right)\in X\times X,$ and $\lambda \in A$, $\left\{d\left(\tilde{x},\tilde{y}\right)\left(\lambda \right):\tilde{x}\left(\lambda \right)=\xi ,\tilde{y}\left(\lambda \right)=\eta \right\}$ is a singleton set, and if for$\ \lambda \in A$, $d_{\lambda }:X\times X\to {{\mathcal R}}^+$ is defined by $d_{\lambda }\left(\tilde{x}\left(\lambda \right),\tilde{y}\left(\lambda \right)\right)=d\left(\tilde{x},\tilde{y}\right)\left(\lambda \right)$, $\tilde{x},\tilde{y}\widetilde{\in }\check{X};$ then $d_{\lambda }$ is a metric on $X.$
\end{theorem} 

\begin{definition} \cite{S3}
Let $(\check{X},d)$ be a soft metric space and $\mathcal{B}$ be a non-null collection of soft elements of $\check{X}$. Then $\mathcal{B}$ is said to be `\textit{open in } $\check{X}$ \textit{ with respect to } $d$' or ` \textit{open in } $(\check{X},d)$' if all elements of $\mathcal{B}$ are interior elements of $\mathcal{B}$. 
\end{definition}

\begin{definition} \cite{S3}
Let $(\check{X},d)$ be a soft metric space and $(Y,A)$ be a non-null soft subset $\in $ ${\mathcal S}\left(\check{X}\right)\ $in $(\check{X},d)$. Then $(Y,A)$ is said to be `\textit{soft open in } $\check{X}$ \textit{ with respect to } $d$' if there is a collection $\mathcal{B}$ of soft elements of $(Y,A)$ such that $\mathcal{B}$ is open with respect to ${\rm d}$ and $\left(Y,A\right)=SS(\mathcal{B})$. 
\end{definition}

\begin{definition} \cite{S3}
 Let $(\check{X},d)$ be a soft metric space. A soft set $(Y,A)\in {\mathcal S}(\check{X})$, is said to be `\textit{soft closed in } $\check{{\check X}}$ \textit{ with respect to } $d$' if its complement ${(Y,A)}^c$ is a member of ${\mathcal S}(\check{X})$ and is soft open in $(\check{X},d)$.
\end{definition} 

\begin{definition} \cite{S3}
Let $(\check{X},d)$ be a soft metric space and $\mathcal{B}$ be a collection of soft elements of $\check{X}$. A soft element $\tilde{a}\widetilde{\in }\mathcal{B}$ is said to be a \textit{limit element }of $\mathcal{B}$, if every open ball $B\left(\tilde{a},\tilde{r}\right)$ containing $\tilde{a}\ $in $(\check{X},d)$ contains at least one element of $\mathcal{B}$ different from $\tilde{a}$.

\noindent The set of all limit elements of $\mathcal{B}$ is said to be the derived set of $\mathcal{B}$ and is denoted by ${\mathcal{B}}^d.$
\end{definition}

\begin{definition} \cite{S3}
Let $(\check{X},d)$ be a soft metric space and $(Y,A)$ $\in {\mathcal S}(\check{X})$. A soft element $\tilde{a}\widetilde{\in }\check{X}$ is said to be a \textit{soft limit element } of $(Y,A)$, if every open ball $B\left(\tilde{a},\tilde{r}\right)$ containing $\tilde{a}$ in $(\check{X},d)$ contains at least one soft element of $(Y,A)$ different from $\tilde{a}$.

\noindent A soft limit element of a soft set $(Y,A)$ may or may not belong to the soft set $(Y,A)$.

\noindent The set of all soft limit elements of $(Y,A)$ is said to be the derived set of $(Y,A)$ and is denoted by ${(Y,A)}^d.$
\end{definition} 

\begin{definition} \cite{S3}
Let $\left(\check{X},d\right)$ be a soft metric space and $\mathcal{B}$ be a collection of soft elements of $\check{X}$. Then the collection of all soft elements of $\mathcal{B}$ and limit elements of $\mathcal{B}$ in $\left(\check{X},d\right)$ is said to be the \textit{closure} of $\mathcal{B}$ in $\left(\check{X},d\right)$. It is denoted by $\widetilde{\mathcal{B}}$.
\end{definition}

\begin{definition} \cite{S3}
Let $\left(\check{X},d\right)$ be a soft metric space and$(Y,A)$ be a soft subset $\in {\mathcal S}(\check{X})$. Then the collection of all soft elements of $(Y,A)$ and soft limit elements of $(Y,A)$ in $\left(\check{X},d\right)$ is said to be the \textit{soft closure} of $(Y,A)$ in $\left(\check{X},d\right)$. It is denoted by $\overline{(Y,A)}$.
\end{definition}

\begin{definition} \cite{S3}
Let $\{{\tilde{x}}_n\}$ be a sequence of soft elements in a soft metric space $\left(\check{X},d\right)$. The sequence $\{{\tilde{x}}_n\}$ is said to be convergent in $\left(\check{X},d\right)$ if there is a soft element $\tilde{x}\widetilde{\in }\check{X}$ such that $d({\tilde{x}}_n,\tilde{x})\to \overline{0}$ as $n\to \infty $.

\noindent This means for every $\widetilde{\varepsilon }\tilde{>}\overline{0}$, chosen arbitrarily, $\exists $ a natural number $N=N(\widetilde{\varepsilon })$, such that $\overline{0}\widetilde{\le }d({\tilde{x}}_n,\tilde{x})\tilde{<}\widetilde{\varepsilon }$ , whenever $n>N$.

\noindent i.e., $n>N\Longrightarrow {\tilde{x}}_n\in B(\tilde{x},\ \widetilde{\varepsilon })$. We denote this by ${\tilde{x}}_n\to \tilde{x}$ as $n\to \infty $ or by ${\mathop{\lim }_{n\to \infty } {\tilde{x}}_n\ }=\tilde{x}$. $\tilde{x}\ $ is said to be the limit of the sequence ${\tilde{x}}_n$ as $n\to \infty $.
\end{definition} 

\section{Soft Vector/Linear Spaces}

\noindent Let $V$ be a vector space over a field $K$ and let $A$ be a parameter set. A soft set $(F,A)$ where $F:A\to \wp (V)$ will be denoted by $F$ only.

\begin{definition}{(Sums and Scalar products of soft sets)}
Let$F_1,F_2,\dots \dots ,F_n$ be $n$ soft sets in $\left(V,A\right).$ Then $F=F_1+F_2+\dots +F_n$ is a soft set over $\left(V,A\right)$ and is defined as $F\left(\lambda \right)=\left\{x_1+x_2+\dots +x_n;x_i\in F_i\left(\lambda \right),i=1,2,\dots ,n\right\},\ \forall \lambda \in A.$

\noindent Let $\alpha \in K$ be a scalar and $F$ be a soft sets over $\left(V,A\right),$ then$\alpha F$ is a soft set over $\left(V,A\right)$ and is defined as follows: $\alpha F=G,\ \ \ G\left(\lambda \right)=\{\alpha x;x\in F(\lambda )\}$, $\lambda \in A.$
\end{definition} 

\begin{definition}
Let $V$ be a vector space over a field $K$ and let $A$ be a parameter set. Let $G$ be a soft set over $\left(V,A\right)$. Now $G$ is said to be a soft vector space or soft linear space of $V$ over $K$ if $G(\lambda )$ is a vector subspace of $V,\ \forall \lambda \in A.$
\end{definition} 

\begin{example} \label{ex424}
Consider the Euclidian n-dimensional space ${{\mathcal R}}^n$ over ${\mathcal R}.$ Let $A=\{1,2,3,\dots ,n\}$ be the set of parameters. Let $G:A\to \wp ({{\mathcal R}}^n)$ be defined as follows:
\[G\left(i\right)=\left\{t\in {{\mathcal R}}^n;i-th\ co-ordinate\ of\ t\ is\ 0\right\},\ i=1,2,\dots ,n.\] 
Then $G$ is a soft vector space or soft linear space of ${{\mathcal R}}^n$ over ${\mathcal R}.$
\end{example}

\begin{proposition}
$\alpha \left(F+G\right)=\alpha F+\alpha G$ for all soft sets $F,G$ over $\left(V,A\right)$ and $\alpha \in K$.
\end{proposition}

\textbf{Proof.}
$\left[\alpha \left(F+G\right)\right]\left(\lambda \right)=\left\{\alpha z;z\in \left(F+G\right)\left(\lambda \right)\right\}$

\noindent $=\left\{\alpha \left(x+y\right);x\in F\left(\lambda \right),y\in G\left(\lambda \right)\right\}$

\noindent $=\{\alpha x+\alpha y;x\in F\left(\lambda \right),y\in G\left(\lambda \right)\}$, $($Since $V$ is a vector space$)$

\noindent Again $\left(\alpha F+\alpha G\right)\left(\lambda \right)=\{x^/+y^/;x^/\in \alpha F\left(\lambda \right),y^/\in \alpha G\left(\lambda \right)\}$

\noindent $=\{\alpha x^{//}+\alpha y^{//};x^{//}\in F\left(\lambda \right),y^{//}\in G\left(\lambda \right)\}.$

\noindent Hence the result follows.

\begin{lemma}
 Let $F_1,F_2,\dots \dots ,F_n$, $G_1,G_2,\dots \dots ,G_m$ be soft sets over $\left(V,A\right)$ and let $F=F_1+F_2+\dots +F_n$ and $G=G_1+G_2+\dots +G_m.$ Then $H=F+G,$ where $H=F_1+F_2+\dots +F_n+G_1+G_2+\dots +G_m.$
\end{lemma} 

\textbf{Proof.}
 The proof is straightforward.

\begin{definition} \label{def401}
Let $x\in V$ and $F$ be a soft setover $\left(V,A\right).$ Then $x+F$ is a soft set over $\left(V,A\right)$ defined as follows: 
\[\left(x+F\right)\left(\lambda \right)=\left\{x+y;y\in F\left(\lambda \right)\right\},\lambda \in A.\] 
\end{definition}

\begin{lemma}
Let $U$ be an ordinary subset of $V$ and let $F$ be a soft set over $\left(V,A\right).$ Then $U+F$ is a soft set over $\left(V,A\right)$ defined as follows:
\[\left(U+F\right)\left(\lambda \right)=\bigcup_{x\in U}{\left\{x+y;y\in F\left(\lambda \right)\right\}},\lambda \in A.\] 
\[\ \ U+F=\bigcup_{x\in U}{(x+F)}\] 
\end{lemma}

\textbf{Proof.}
Follows from Definition \ref{def401} and the Definition of union of soft sets.

\begin{definition}{(Soft Vector Subspaces)}
Let $F$ be a soft vector space of $V$ over $K$. Let $G:A\to \wp (V)$ be a soft set over $\left(V,A\right)$. Then $G$ is said to be a soft vector subspace of $F$ if 

\begin{itemize}
\item [{(i).}] for each $\lambda \in A,G(\lambda )$ is a vector subspace of $V$ over $K$ and 

\item [{(ii).}] $F\left(\lambda \right)\supseteq G\left(\lambda \right),\ \forall \lambda \in A.$
\end{itemize}
\end{definition} 

\begin{theorem}
 A soft subset $G$ of a soft vector space$F$ is a soft vector subspace of $F$ if and only if for all scalars $\alpha $,$\beta \in K,\alpha G+\beta G\subset G.$
\end{theorem}

\textbf{Proof.}
 Let $G$ be a soft vector subspace of $F$ of $V$ over $K$.

\noindent Let $\lambda \in A$, $\left(\alpha G+\beta G\right)\left(\lambda \right)=\{x^/+y^/;x^/\in \alpha G\left(\lambda \right),y^/\in \beta G\left(\lambda \right)\}$

\noindent $=\{\alpha x+\beta y;x,y\in G(\lambda )\}\subset G(\lambda ),$

\noindent [Now, $G(\lambda )$ is a vector subspace over $K$ and $x,y\in G(\lambda )$,$\ \alpha $,$\beta \in K\Rightarrow \alpha x+\beta y\in G(\lambda )$]

\noindent $\ \ \alpha G+\beta G\subset G,$ the given condition is satisfied.

\noindent 

\noindent Conversely, let the given condition hold.

\noindent We have, $\left(\alpha G+\beta G\right)\left(\lambda \right)=\{\alpha x+\beta y;x,y\in G(\lambda )\}\ ,\ \forall \lambda \in A.$

\noindent By the given condition, $\alpha G+\beta G\subset G$ i.e., $\{\alpha x+\beta y;x,y\in G(\lambda )\}\subset G(\lambda ),\forall \lambda \in A.$

\noindent i.e., for all $x,y\in G\left(\lambda \right)$ and for every scalar $\alpha $,$\beta \in K$, $\alpha x+\beta y\in G\left(\lambda \right)$

\noindent $\Rightarrow G(\lambda )$ is a vector subspace of $F$ over $K.$ This is true for all $\lambda \in A.$

\noindent Also since $G$ is a soft subset of $F,\ F\left(\lambda \right)\supseteq G\left(\lambda \right),\ \forall \lambda \in A.$

\noindent Hence $G$ is a soft vector subspace of $F$.

\begin{proposition}
 If $F$ and $G$ be two soft vector subspaces of $H$ over $K$ and $\alpha $ be a scalar, then $F+G$ and $\alpha F$ are soft vector subspaces of $H$ over $K$. 
\end{proposition} 

\textbf{Proof.}
 The proof is straight forward.

\begin{proposition}
 If $\{F_i\}$ be a family of soft vector subspaces of $H$ over $K$, then $G=\bigcap_i{F_i}$ is a soft vector subspaces of $H$ over $K$. 
\end{proposition} 

\textbf{Proof.}
 The proof is straight forward.

\section{Soft vectors in soft vector spaces}

\noindent In this section we introduce the concept of soft vectors in soft vector spaces and study some of their basic properties.

\begin{definition}
 Let $G$ be a soft vector space of $V$ over $K$. Then a soft element of $V$ is said to be a soft vector of $G.$In a similar manner a soft element of the soft set $(K,A)$ is said to be a soft scalar, K being the scalar field.
\end{definition}

\begin{example}
Consider the soft vector space $G$ as Example \ref{ex424}. Let $\tilde{x}$ be a soft element of $G$ as the following; 

\noindent $\tilde{x}\left(i\right)=\left(1,1,..,0_{i-th},..,1\right)\in {{\mathcal R}}^n,i=1,2,..,n.$ Then $\tilde{x}$ is a soft vector of $G.$
\end{example}

\begin{definition}
 A soft vector $\tilde{x}$ in a soft vector space $G$ is said to be the null soft vector if $\tilde{x}\left(\lambda \right)=\theta ,\ \forall \lambda \in A$, $\theta $ being the zero element of $V$. It will be denoted by $\Theta $. A soft vector is said to be non-null if itis not a null soft vector.
\end{definition} 

\begin{definition}
Let $\tilde{x},\tilde{y}$ be soft vectors of $G$ and $\tilde{k}$ be a soft scalar. Then the addition $\tilde{x}+\tilde{y}$ of $\tilde{x},\tilde{y}$ and scalar multiplication $\tilde{k}.\ \tilde{x}$ of $\tilde{k}$ and $\tilde{x}$ are defined by 

\noindent $\left(\tilde{x}+\tilde{y}\right)\left(\lambda \right)=\tilde{x}\left(\lambda \right)+\tilde{y}\left(\lambda \right),\ \ \left(\tilde{k}.\ \tilde{x}\right)\left(\lambda \right)=\tilde{k}\left(\lambda \right).\tilde{x}\left(\lambda \right),\ \forall \lambda \in A.$ Obviously, $\tilde{x}+\tilde{y},\tilde{k}.\ \tilde{x}$ are soft vectors of $G.$
\end{definition}

\begin{theorem} 
In a soft vector space $G$ of $V$ over $K$,

\begin{itemize}
\item [{(i).}] $\overline{0}.\widetilde{\alpha }=\Theta ,\ $ for all $\widetilde{\alpha }\widetilde{\in }G;$

\item [{(ii).}]  $\tilde{k}.\Theta =\Theta ,\ $ for all soft scalar $\tilde{k}.$

\item [{(iii).}]  $\left(-\overline{1}\right).\ \widetilde{\alpha }=-\widetilde{\alpha }$$,$ for all $\widetilde{\alpha }\widetilde{\in }\ G.$

\end{itemize}
\end{theorem}

\textbf{Proof.}
(i) We have, $\left(\overline{0}.\widetilde{\alpha }\right)\left(\lambda \right)=\overline{0}\left(\lambda \right).\widetilde{\alpha }\left(\lambda \right)=0.\widetilde{\alpha }\left(\lambda \right)=\ \theta ,\ \ \forall \lambda \in A.$

\noindent $\Rightarrow \overline{0}.\widetilde{\alpha }=\Theta ,\ $ for all $\widetilde{\alpha }\widetilde{\in }\ G.$

\noindent (ii) $\left(\tilde{k}.\Theta \right)\left(\lambda \right)=\tilde{k}\left(\lambda \right).\Theta \left(\lambda \right)=\tilde{k}\left(\lambda \right).\theta =\theta ,\ \ \forall \lambda \in A.$

\noindent $\Rightarrow \tilde{k}.\Theta =\Theta ,\ $ for all soft scalar $\tilde{k}.$

\noindent (iii) $\left(\left(-\overline{1}\right).\ \widetilde{\alpha }\right)\left(\lambda \right)=\left(-\overline{1}\right)\left(\lambda \right).\widetilde{\alpha }\left(\lambda \right)=\left(-1\right).\widetilde{\alpha }\left(\lambda \right)=-\widetilde{\alpha }\left(\lambda \right)=(-\widetilde{\alpha })\left(\lambda \right),\ $ 

\noindent $\forall \lambda \in A.$

\noindent $\Rightarrow \left(-\overline{1}\right).\ \widetilde{\alpha }=-\widetilde{\alpha },$ for all $\widetilde{\alpha }\widetilde{\in }\ G.$

\begin{remark}
 However, $\tilde{k}.\widetilde{\alpha }=\Theta $ does not necessarily imply that either $\tilde{k}=\overline{0}$ or $\widetilde{\alpha }=\Theta $. For example let us consider the soft vector space as Example \ref{ex424} Let $\tilde{k}\left(1\right)=1,$ and $\tilde{k}\left(i\right)=0,$ for $i=2,3,..,n\ $ and $\widetilde{\alpha }\left(1\right)=\theta ,$ and $\widetilde{\alpha }\left(i\right)=\left(1,1,..,0_{i-th},..,1\right)\in {{\mathcal R}}^n,i=2,3,..,n.$ Then $\left(\tilde{k}.\widetilde{\alpha }\right)\left(1\right)=\tilde{k}\left(1\right).\ \widetilde{\alpha }\left(1\right)=1.\ \theta =\theta =\Theta \left(1\right)$ and 

\noindent $\left(\tilde{k}.\widetilde{\alpha }\right)\left(i\right)=\tilde{k}\left(i\right).\ \widetilde{\alpha }\left(i\right)=0.\ \left(1,1,..,0_{i-th},..,1\right)=\theta =\Theta \left(i\right),$ for $i=2,3,..,n$.

\noindent $\ \tilde{k}.\widetilde{\alpha }=\Theta ,$ but neither $\tilde{k}=\overline{0}$ nor $\widetilde{\alpha }=\Theta $.
\end{remark} 

\begin{theorem}
 A non-null soft subset $(W,A)$ of a soft vector space $G$ of $V$ over $K$, is a soft subspace of $G$ if and only if $\widetilde{\alpha },\ \widetilde{\beta }\widetilde{\in }(W,A)$ and $\tilde{k},\ \tilde{s}$ be soft scalars then  $\tilde{k}.\widetilde{\alpha }+\tilde{s}.\ \widetilde{\beta }\widetilde{\in }(W,A)$.
\end{theorem}

\textbf{Proof.}
 Let $(W,A)$ be a soft vector sub space of $G$ of $V$ over $K$. Let $\widetilde{\alpha },\ \widetilde{\beta }\widetilde{\in }(W,A)$ and $\tilde{k},\ \tilde{s}$ be soft scalars, then 

\noindent $\left(\tilde{k}.\widetilde{\alpha }+\tilde{s}.\ \widetilde{\beta }\right)\left(\lambda \right)=\tilde{k}\left(\lambda \right).\widetilde{\alpha }\left(\lambda \right)+\tilde{s}\left(\lambda \right).\widetilde{\beta }\left(\lambda \right)\in W\left(\lambda \right),\ \forall \lambda \in A.$ ($W\left(\lambda \right)$ is a vector subspace of $V$ for each $\lambda \in A,\ \tilde{k}\left(\lambda \right),\ \tilde{s}\left(\lambda \right)\in K$, $\widetilde{\alpha }\left(\lambda \right)+\widetilde{\beta }\left(\lambda \right)\in W\left(\lambda \right),\ \forall \lambda \in A.$)
\[\ \tilde{k}.\widetilde{\alpha }+\tilde{s}.\ \widetilde{\beta }\widetilde{\in }\left(W,A\right).\]

\noindent Conversely, let the given condition be satisfied.

\noindent Then for all soft scalars $\tilde{k},\ \tilde{s}$ and soft vectors $\widetilde{\alpha },\ \widetilde{\beta }\widetilde{\in }\left(W,A\right),\ \tilde{k}.\widetilde{\alpha }+\tilde{s}.\ \widetilde{\beta }\widetilde{\in }\left(W,A\right).$

\noindent i.e., $\left(\tilde{k}.\widetilde{\alpha }+\tilde{s}.\ \widetilde{\beta }\right)\left(\lambda \right)\in W\left(\lambda \right),\ \forall \lambda \in A$, i.e., $\tilde{k}\left(\lambda \right).\widetilde{\alpha }\left(\lambda \right)+\tilde{s}\left(\lambda \right).\widetilde{\beta }\left(\lambda \right)\in W\left(\lambda \right),\ $

\noindent $\forall \lambda \in A$;

\noindent $\Rightarrow \ W\left(\lambda \right)$ is a vector subspace of $V$ for each $\lambda \in A.$ Also it is obvious that, $W\left(\lambda \right)\subset V\left(\lambda \right),\ \forall \lambda \in A.$

\noindent Hence, $(W,A)$ is a soft vector sub space of $G$ of $V$ over $K$.

\begin{definition}
Let $G$ be a soft vector space of $V$ over $K$. Let ${\widetilde{\alpha }}_1,{\widetilde{\alpha }}_2,\dots ,$

\noindent ${\widetilde{\alpha }}_n\widetilde{\in }G.$ A soft vector $\widetilde{\beta }$ in $G$ is said to be a linear combination of the soft vectors ${\widetilde{\alpha }}_1,{\widetilde{\alpha }}_2,\dots ,{\widetilde{\alpha }}_n$ if $\widetilde{\beta }$ can be expressed as $\widetilde{\beta }={\tilde{c}}_1.\ {\widetilde{\alpha }}_1+{\tilde{c}}_2.\ {\widetilde{\alpha }}_2+..+{\tilde{c}}_n.\ {\widetilde{\alpha }}_n,$ for some soft scalars ${\tilde{c}}_1,\ {\tilde{c}}_2,\ \dots ,{\tilde{c}}_n$.
\end{definition} 

\begin{example}
Consider the soft vector space $G$ as Example \ref{ex424}. 

\noindent Let ${\widetilde{\alpha }}_i=\left(1,1,..,0_{i-th},..,1\right)\in {{\mathcal R}}^n,i=1,2,3.$ Then ${\widetilde{\alpha }}_1+\ {\widetilde{\alpha }}_2,\ {\widetilde{\alpha }}_1+\ {\widetilde{\alpha }}_2+{\widetilde{\alpha }}_3,\ {\overline{2}.\widetilde{\alpha }}_1+\overline{5}.\ {\widetilde{\alpha }}_2+{\widetilde{\alpha }}_3$ are linear combinations of ${\widetilde{\alpha }}_1,\ {\widetilde{\alpha }}_2,{\widetilde{\alpha }}_3.$
\end{example} 

\begin{definition}
 A finite set of soft vectors $\{{\widetilde{\alpha }}_1,{\widetilde{\alpha }}_2,\dots ,{\widetilde{\alpha }}_n\}$ of a soft vector space $G$ is said to be linearly dependent in $G$ if there exists soft scalars ${\tilde{c}}_1,\ {\tilde{c}}_2,\ \dots ,{\tilde{c}}_n$ not all $\overline{0}$ such that

\begin{equation} \label{eq426}
{\tilde{c}}_1.\ {\widetilde{\alpha }}_1+{\tilde{c}}_2.\ {\widetilde{\alpha }}_2+..+{\tilde{c}}_n.\ {\widetilde{\alpha }}_n=\Theta 
\end{equation}

\noindent The set is said to be linearly independent in $G$ if the equality \ref{eq426} is satisfied only when ${\tilde{c}}_1=\ {\tilde{c}}_2=\ \dots ={\tilde{c}}_n=\overline{0}$.

\noindent An arbitrary set $S$ of soft vectors of $G$ is said to be linearly dependent in $G$ if there exists a finite subset of $S$ which is linearly dependent in $G$.

\noindent A set of soft vectors which is not linearly dependent is said to be linearly independent.
\end{definition} 

\begin{proposition}
A set $S=\{{\widetilde{\alpha }}_1,{\widetilde{\alpha }}_2,\dots ,{\widetilde{\alpha }}_n\}$ of soft vectors in a soft vector space $G$ over $V$ is linearly independent if and only if the sets 

\noindent $S\left(\lambda \right)=\{{\widetilde{\alpha }}_1\left(\lambda \right),{\widetilde{\alpha }}_2\left(\lambda \right),\dots ,{\widetilde{\alpha }}_n\left(\lambda \right)\}$ are linearly independent in $V,\forall \lambda \in A.$
\end{proposition} 

\textbf{Proof.} 
Let $S$ be linearly independent. Then for any set of soft scalars ${\tilde{c}}_1,\ {\tilde{c}}_2,\ \dots ,{\tilde{c}}_n,$

\begin{equation} \label{eq419}
{\tilde{c}}_1.\ {\widetilde{\alpha }}_1+{\tilde{c}}_2.\ {\widetilde{\alpha }}_2+..+{\tilde{c}}_n.\ {\widetilde{\alpha }}_n=\Theta \Leftrightarrow {\tilde{c}}_1=\ {\tilde{c}}_2=\ \dots ={\tilde{c}}_n=\overline{0}.
\end{equation}

\noindent Let ${\lambda }_0\in A,$ and $S\left({\lambda }_0\right)=\{{\widetilde{\alpha }}_1\left({\lambda }_0\right),{\widetilde{\alpha }}_2\left({\lambda }_0\right),\dots ,{\widetilde{\alpha }}_n\left({\lambda }_0\right)\}.$ Let $c_1.\ {\widetilde{\alpha }}_1\left({\lambda }_0\right)+c_2.\ {\widetilde{\alpha }}_2\left({\lambda }_0\right)+..+c_n.\ {\widetilde{\alpha }}_n\left({\lambda }_0\right)=\theta $. Let us consider any set of soft scalars ${\tilde{c}}_1,\ {\tilde{c}}_2,\ \dots ,{\tilde{c}}_n,$ such that ${\tilde{c}}_1\left({\lambda }_0\right)=c_1,{\tilde{c}}_2\left({\lambda }_0\right){=c}_2..{\tilde{c}}_n\left({\lambda }_0\right)=c_n;$ then since $S$ is linearly independent, from (\ref{eq419}), it follows that,${\tilde{c}}_1=\ {\tilde{c}}_2=\ \dots ={\tilde{c}}_n=\overline{0}$. Hence ${\tilde{c}}_1\left({\lambda }_0\right)={\tilde{c}}_2\left({\lambda }_0\right)=\dots ={\tilde{c}}_n\left({\lambda }_0\right)=0$ i.e., $c_1=c_2=\ \dots =c_n=0$. Hence $S\left({\lambda }_0\right)=\{{\widetilde{\alpha }}_1\left({\lambda }_0\right),{\widetilde{\alpha }}_2\left({\lambda }_0\right),\dots ,{\widetilde{\alpha }}_n\left({\lambda }_0\right)\}$ is linearly independent in $V.$

\noindent Since ${\lambda }_0\in A,$ is arbitrary, it follows that $S\left(\lambda \right)=\{{\widetilde{\alpha }}_1\left(\lambda \right),{\widetilde{\alpha }}_2\left(\lambda \right),\dots ,{\widetilde{\alpha }}_n\left(\lambda \right)\},$ are linearly independent in $V,\forall \lambda \in A.$

\noindent 

\noindent Conversely let $S\left(\lambda \right)=\{{\widetilde{\alpha }}_1\left(\lambda \right),{\widetilde{\alpha }}_2\left(\lambda \right),\dots ,{\widetilde{\alpha }}_n\left(\lambda \right)\},\ $ be linearly independent in $V,\forall \lambda \in A.$

\noindent Let ${\tilde{c}}_1,\ {\tilde{c}}_2,\ \dots ,{\tilde{c}}_n$ be any set of soft scalars such that ${\tilde{c}}_1.\ {\widetilde{\alpha }}_1+{\tilde{c}}_2.\ {\widetilde{\alpha }}_2+..+{\tilde{c}}_n.\ {\widetilde{\alpha }}_n=\Theta $.Then ${\tilde{c}}_1\left(\lambda \right).\ {\widetilde{\alpha }}_1\left(\lambda \right)+{\tilde{c}}_2\left(\lambda \right).\ {\widetilde{\alpha }}_2\left(\lambda \right)+..+{\tilde{c}}_n\left(\lambda \right).\ {\widetilde{\alpha }}_n\left(\lambda \right)=\theta \Rightarrow {\tilde{c}}_1\left(\lambda \right)=\ {\tilde{c}}_2\left(\lambda \right)=\ \dots ={\tilde{c}}_n\left(\lambda \right)=0,\ \forall \lambda \in A\Rightarrow {\tilde{c}}_1=\ {\tilde{c}}_2=\ \dots ={\tilde{c}}_n=\overline{0}$.

\noindent Hence $S=\{{\widetilde{\alpha }}_1,{\widetilde{\alpha }}_2,\dots ,{\widetilde{\alpha }}_n\}$ is linearly independent.

\begin{proposition}
A set $S=\{{\widetilde{\alpha }}_1,{\widetilde{\alpha }}_2,\dots ,{\widetilde{\alpha }}_n\}$ of soft vectors in a soft vector space $G$ over $V$ is linearly dependent if and only if the sets 

\noindent $S\left(\lambda \right)=\{{\widetilde{\alpha }}_1\left(\lambda \right),{\widetilde{\alpha }}_2\left(\lambda \right),\dots ,{\widetilde{\alpha }}_n\left(\lambda \right)\},\ $ are linearly dependent in $V$ for some $\lambda \in A$.
\end{proposition} 

\textbf{Proof.}
 Let $S$ be linearly dependent. Then there is a set of soft scalars ${\tilde{c}}_1,\ {\tilde{c}}_2,\ \dots ,{\tilde{c}}_n,$ not all equal to $\overline{0}$ such that ${\tilde{c}}_1.\ {\widetilde{\alpha }}_1+{\tilde{c}}_2.\ {\widetilde{\alpha }}_2+..+{\tilde{c}}_n.\ {\widetilde{\alpha }}_n=\Theta $. Then ${\tilde{c}}_1\left(\lambda \right).\ {\widetilde{\alpha }}_1\left(\lambda \right)+{\tilde{c}}_2\left(\lambda \right).\ {\widetilde{\alpha }}_2\left(\lambda \right)+..+{\tilde{c}}_n\left(\lambda \right).\ {\widetilde{\alpha }}_n\left(\lambda \right)=\theta $, $\forall \lambda \in A$ and there is at least one ${\lambda }_0\in A,$ such that ${\tilde{c}}_1\left({\lambda }_0\right),\ {\tilde{c}}_2\left({\lambda }_0\right),\dots .,{\tilde{c}}_n\left({\lambda }_0\right)$ are not all zeros. Then ${\tilde{c}}_1\left({\lambda }_0\right).\ {\widetilde{\alpha }}_1\left({\lambda }_0\right)+{\tilde{c}}_2\left({\lambda }_0\right).\ {\widetilde{\alpha }}_2\left({\lambda }_0\right)+..+{\tilde{c}}_n\left({\lambda }_0\right).\ {\widetilde{\alpha }}_n\left({\lambda }_0\right)=\theta $ and ${\tilde{c}}_1\left({\lambda }_0\right),\ {\tilde{c}}_2\left({\lambda }_0\right),\dots .,{\tilde{c}}_n\left({\lambda }_0\right)$ are not all zeros. Proving that $S\left({\lambda }_0\right)=\{{\widetilde{\alpha }}_1\left({\lambda }_0\right),{\widetilde{\alpha }}_2\left({\lambda }_0\right),\dots ,{\widetilde{\alpha }}_n\left({\lambda }_0\right)\}$ is linearly dependent.

\noindent 

\noindent Conversely let $S\left({\lambda }_0\right)=\{{\widetilde{\alpha }}_1\left({\lambda }_0\right),{\widetilde{\alpha }}_2\left({\lambda }_0\right),\dots ,{\widetilde{\alpha }}_n\left({\lambda }_0\right)\},\ $ be linearly dependent for some ${\lambda }_0\in A$.. Then there is a set of scalars $c_1,c_2,\ \dots ,c_n$ not all zeros such that $c_1.\ {\widetilde{\alpha }}_1\left({\lambda }_0\right)+c_2.\ {\widetilde{\alpha }}_2\left({\lambda }_0\right)+..+c_n.\ {\widetilde{\alpha }}_n\left({\lambda }_0\right)=\theta .$

\noindent Let ${\tilde{c}}_1,\ {\tilde{c}}_2,\ \dots ,{\tilde{c}}_n$ be a set of soft scalars such that${\tilde{c}}_i\left({\lambda }_0\right)=c_i,$ and ${\tilde{c}}_i\left(\lambda \right)=0$for $\lambda \in A\backslash \{{\lambda }_0\}$, for $i=1,2,\dots ,n.$ Then ${\tilde{c}}_1,\ {\tilde{c}}_2,\ \dots ,{\tilde{c}}_n,$ are not all equal to $\overline{0}$ and ${\tilde{c}}_1.\ {\widetilde{\alpha }}_1+{\tilde{c}}_2.\ {\widetilde{\alpha }}_2+..+{\tilde{c}}_n.\ {\widetilde{\alpha }}_n=\Theta $. Hence $S$ is linearly dependent.

\section{Soft norm and soft normed linear spaces}

\noindent Let $X$ be a vector space over a field $K$,$X$ is also our initial universe set and $A$ be a non-empty set of parameters. Let $\check{X}$ be the absolute soft vector space i.e., $\check{X}\left(\lambda \right)=X$, $\forall \lambda \in A$. We use the notation $\tilde{x},\ \tilde{y},\ \tilde{z}$ to denote soft vectors of a soft vector space and $\tilde{r},\ \tilde{s},\ \tilde{t}$ to denote soft real numbers whereas $\overline{r},\ \overline{s},\ \overline{t}$ will denote a particular type of soft real numbers such that $\overline{r}\left(\lambda \right)=r,\ $ for all $\lambda \in A$ etc. For example $\overline{0}$ is the soft real number such that $\overline{0}\left(\lambda \right)=0,\ $ for all $\lambda \in A$. Note that, in general, $\tilde{r}$ is not related to $r$.

\subsection{Definitions and examples of soft norm and soft normed linear spaces}

\begin{definition}
 Let $\check{X}$ be the absolute soft vector space i.e., $\check{X}\left(\lambda \right)=X$, $\forall \lambda \in A$. Then a mapping $\left\|.\right\|:SE(\check{X})\to {R\left(A\right)}^*$  is said to be a soft norm on the soft vector space $\check{X}$ if $\left\|.\right\|$ satisfies the following conditions:

\begin{itemize}

\item [{(N1).}] $\left\|\tilde{x}\right\|\widetilde{\ge }\overline{0},\ $ for all $\tilde{x}\widetilde{\in }\check{X}$;

\item [{(N2).}] $\left\|\tilde{x}\right\|=\overline{0}$ if and only if $\tilde{x}=\Theta ;$

\item [{(N3).}]  $\left\|\widetilde{\alpha }.\tilde{x}\right\|=|\widetilde{\alpha }|\left\|\tilde{x}\right\|$ for all $\tilde{x}\widetilde{\in }\check{X}$ and for every soft scalar $\widetilde{\alpha }$;

\item [{(N4).}] For all $\tilde{x},\tilde{y}\widetilde{\in }\check{X},\left\|\tilde{x}+\tilde{y}\right\|\widetilde{\le }\left\|\tilde{x}\right\|+\left\|\tilde{y}\right\|$.
\end{itemize}

\noindent The soft vector space $\check{X}$ with a soft norm $\left\|.\right\|$ on $\check{X}$ is said to be a soft normed linear space and is denoted by $(\check{X},\left\|.\right\|,A)$ or $(\check{X},\left\|.\right\|)$. (N1), (N2), (N3) and (N4) are said to be soft norm axioms.

\end{definition}

\begin{example} \label{ex411}
Let ${\mathcal R}(A)$ be the set of all soft real numbers. We define $\left\|.\right\|:{\mathcal R}(A)\to {{\mathcal R}\left(A\right)}^*$ ,by,

\noindent $\left\|\tilde{x}\right\|=|\tilde{x}|,\ $ for all $\tilde{x}\widetilde{\in }{\rm \ }{\mathcal R}{\rm (}A)$, where $|\tilde{x}|$ denotes the modulus of soft real numbers. Then $\left\|.\right\|$ satisfied all the soft norm axioms so, $\left\|.\right\|$ is a soft norm on ${\mathcal R}(A)$ and $({\mathcal R}(A),\left\|.\right\|,A)$ or${\rm (}{\mathcal R}{\rm (}A),\left\|.\right\|)$ is a soft normed linear space.
\end{example} 

\begin{example} \label{ex421}
Every parametrized family of crisp norms\{${\left\|.\right\|}_{\lambda }:\lambda \in A$\} on a crisp vector space $X$ can be considered as a soft norm on the soft vector space $\check{X}$.
\end{example}

\textbf{Proof.}
Let $\check{X}$ be the absolute soft vector space over a field $K$, $A$ be a non-empty set of parameters. Let \{${\left\|.\right\|}_{\lambda }:\lambda \in A$\} be a family of crisp norms on the vector space $X$. Let $\tilde{x}\widetilde{\in }\check{X}$, then  $\tilde{x}(\lambda )\in X$, for every $\lambda \in A$. Let us define a mapping $\left\|.\right\|:\check{X}\to {R\left(A\right)}^*$by $\left\|\tilde{x}\right\|\left(\lambda \right)={\left\|\tilde{x}\left(\lambda \right)\right\|}_{\lambda },\ \forall \lambda \in A,\ \forall \tilde{x}\widetilde{\in }\check{X}.$

\noindent Then $\left\|.\right\|$ is a soft norm on $\check{X}$.

\noindent 

\noindent To verify it we now verify the conditions (N1), (N2), (N3) and (N4) for soft norm.

\noindent (N1). We have $\left\|\tilde{x}\right\|\left(\lambda \right)={\left\|\tilde{x}\left(\lambda \right)\right\|}_{\lambda }\ge 0,\forall \lambda \in A,\ \forall \tilde{x}\widetilde{\in }\check{X},$

\noindent $\ \left\|\tilde{x}\right\|\widetilde{\ge }\overline{0},\ $ for all $\tilde{x}\widetilde{\in }\check{X}$.

\noindent (N2). $\left\|\tilde{x}\right\|=\Theta $

\noindent $\Longleftrightarrow \left\|\tilde{x}\right\|\left(\lambda \right)=\theta ,\ \forall \lambda \in A$ 

\noindent $\Longleftrightarrow {\left\|\tilde{x}\left(\lambda \right)\right\|}_{\lambda }=\theta ,\ \forall \lambda \in A$ 

\noindent $\Longleftrightarrow \tilde{x}\left(\lambda \right)=\theta ,\ \forall \lambda \in A$ 

\noindent $\Longleftrightarrow \tilde{x}=\Theta $ 

\noindent Therefore (N2) is satisfied.

\noindent (N3). We have,  $\left\|\alpha .\tilde{x}\right\|\left(\lambda \right)={\left\|\alpha .\tilde{x}\left(\lambda \right)\right\|}_{\lambda }$

\noindent $=\left|\alpha \right|{\left\|\tilde{x}\left(\lambda \right)\right\|}_{\lambda },\ [$ because ${\left\|\alpha .\tilde{x}\left(\lambda \right)\right\|}_{\lambda }=|\alpha |{\left\|\tilde{x}\left(\lambda \right)\right\|}_{\lambda },\ \forall \lambda \in A]$ 

\noindent $=(\left|\alpha \right|\left\|\tilde{x}\right\|)\left(\lambda \right),\ \forall \lambda \in A.$ 

\noindent $\ \left\|\alpha .\tilde{x}\right\|=\left|\alpha \right|\left\|\tilde{x}\right\|$, for all $\tilde{x}\widetilde{\in }\check{X}$ and for every scalar $\alpha \in K$.

\noindent (N4). For all $\tilde{x},\tilde{y}\widetilde{\in }\check{X},$

\noindent $\left[\left\|\tilde{x}\right\|+\left\|\tilde{y}\right\|\right]\left(\lambda \right)=\left\|\tilde{x}\right\|\left(\lambda \right)+\left\|\tilde{y}\right\|\left(\lambda \right)$
 
\noindent $={\left\|\tilde{x}\left(\lambda \right)\right\|}_{\lambda }+{\left\|\tilde{y}\left(\lambda \right)\right\|}_{\lambda }$ 

\noindent  $\ge {\left\|\tilde{x}\left(\lambda \right)+\tilde{y}\left(\lambda \right)\right\|}_{\lambda }$, [by the property of triangle inequality of ${\left\|.\right\|}_{\lambda }$]

\noindent $=\ \left\|\tilde{x}+\tilde{y}\right\|\left(\lambda \right), \forall \lambda \in A.$ 

\noindent Therefore $\ \left\|\tilde{x}\right\|+\left\|\tilde{y}\right\|\widetilde{\ge }\left\|\tilde{x}+\tilde{y}\right\|.$ 

\noindent Thus (N4) is satisfied.

\noindent $\ \left\|.\right\|$ is a soft norm on $\check{X}$ and consequently $(\check{X},\left\|.\right\|)$) is a soft normed linear space.

\begin{proposition}
Every crisp norm ${\left\|.\right\|}_X$ on a crisp vector space $X$ can be extended to a soft norm on the soft vector space $\check{X}$.
\end{proposition}

\textbf{Proof.}
First we construct the absolute soft vector space $\check{X}$ using a non-empty set of parameters $A$.

\noindent Let us define a mapping $\left\|.\right\|:SE(\check{X})\to {R(A)}^*$ by $\left\|\tilde{x}\right\|\left(\lambda \right)={\left\|\tilde{x}\left(\lambda \right)\right\|}_X,\forall \lambda \in A,\ \forall \tilde{x}\widetilde{\in }\check{X}.$

\noindent Then using the same procedure as Example \ref{ex421}, it can be easily proved that $\left\|.\right\|$ is a soft norm on $\check{X}$.

\noindent This soft norm is generated using the crisp norm ${\left\|.\right\|}_X$ and it is said to be the soft norm generated by ${\left\|.\right\|}_X$.

\begin{theorem}{(Decomposition Theorem)}
If a soft norm $\left\|.\right\|$ satisfies the condition

\noindent (N5). For $\xi \in X,$ and $\lambda \in A$, $\left\{\left\|\tilde{x}\right\|\left(\lambda \right):\tilde{x}\left(\lambda \right)=\xi \right\}$ is a singleton set. 

\noindent And if for each $\lambda \in A$, ${\left\|.\right\|}_{\lambda }:X\to R^+$ be a mapping such that for each $\xi \in X,\ {\left\|\xi \right\|}_{\lambda }=\left\|\tilde{x}\right\|\left(\lambda \right)$, where  $\tilde{x}\ \widetilde{\in }\ \check{X}$ such that $\tilde{x}\left(\lambda \right)=\xi $. Then for each $\lambda \in A,$ ${\left\|.\right\|}_{\lambda }$ is a norm on $X.$
\end{theorem} 

\textbf{Proof.}
Clearly ${\left\|.\right\|}_{\lambda }:X\to R^+$ is a rule that assign a vector of $X$ to a non-negative crisp real number $\forall \lambda \in A.$ Now the well defined property of ${\left\|.\right\|}_{\lambda }$ ,$\ \forall \lambda \in A$ is follows from the condition (N5) and the soft norm axioms gives the norm conditions of ${\left\|.\right\|}_{\lambda }$ ,$\ \forall \lambda \in A$. Thus the soft norm satisfying (N5) gives a parameterized family of crisp norms. With this point of view, it also follows that, a soft norm, satisfying (N5) is a particular `soft mapping' as defined by P. Majumdar, et al. in \cite{Pin1} where $\left\|.\right\|:A\to {{{\rm (}{\mathcal R}}^+)}^X$.

\begin{proposition}
Let $(\check{X},\left\|.\right\|,A)$ be a soft normed linear space. Let us define $d:\check{X}\times \check{X}\to {R\left(A\right)}^*$ by $d\left(\tilde{x},\tilde{y}\right)=\left\|\tilde{x}-\tilde{y}\right\|$, for all $\tilde{x},\tilde{y}\widetilde{\in }\check{X}$. Then $d$ is a soft metric on $\check{X}$. 
\end{proposition} 

\textbf{Proof.}
 We have, (M1). $\ d\left(\tilde{x},\tilde{y}\right)=\left\|\tilde{x}-\tilde{y}\right\|\widetilde{\ge }\overline{0},\ $ for all $\tilde{x},\tilde{y}\widetilde{\in }\check{X}$. [using (N1)]

\noindent (M2). $d\left(\tilde{x},\tilde{y}\right)=\ \overline{0}\Longleftrightarrow \left\|\tilde{x}-\tilde{y}\right\|=\ \overline{0}\Longleftrightarrow \tilde{x}=\tilde{y}.$ [using (N2)]

\noindent (M3).$\ d\left(\tilde{x},\tilde{y}\right)=\left\|\tilde{x}-\tilde{y}\right\|=\left\|\tilde{y}-\tilde{x}\right\|=d\left(\tilde{y},\tilde{x}\right)$ for all $\tilde{x},\tilde{y}\widetilde{\in }\check{X}$. [using (N3)]

\noindent (M4). $d\left(\tilde{x},\tilde{y}\right)+d\left(\tilde{y},\tilde{z}\right)=\left\|\tilde{x}-\tilde{y}\right\|+\left\|\tilde{y}-\tilde{z}\right\|\widetilde{\ge }\left\|\tilde{x}-\tilde{z}\right\|=\ d\left(\tilde{x},\tilde{z}\right).$ [using (N4)]

\noindent $\ \ d$ is a soft metric on $\check{X}$.$d$ is said to be the soft metric induced by the soft norm $\left\|.\right\|$. From the above proposition it also follows that every soft normed linear space is also a soft metric space.

\begin{proposition}{(Translation invariance)}
 A soft metric $d$ induced by a soft norm $\left\|.\right\|$ on a normed linear space $(\check{X},\left\|.\right\|)$ satisfies

\begin{itemize} 
\item [{(1).}] $d\left(\tilde{x}+\tilde{a},\ \tilde{y}+\tilde{a}\right)=d(\tilde{x},\tilde{y});$ 

\item [{(2).}] $d\left(\widetilde{\alpha }.\tilde{x},\ \widetilde{\alpha }.\tilde{y}\right)=\left|\widetilde{\alpha }\right|d(\tilde{x},\tilde{y})$, for all $\tilde{x},\ \tilde{y}\widetilde{\in }\check{X}$ and for every soft scalar $\widetilde{\alpha }$.
\end{itemize} 
\end{proposition}

\textbf{Proof.}
 We have, $d\left(\tilde{x}+\tilde{a},\ \tilde{y}+\tilde{a}\right)=\left\|(\tilde{x}+\tilde{a})-(\tilde{y}+\tilde{a})\right\|=\left\|\tilde{x}-\tilde{y}\right\|=\ d(\tilde{x},\tilde{y})$ and 
\[d\left(\widetilde{\alpha }.\tilde{x},\ \widetilde{\alpha }.\tilde{y}\right)=\left\|\widetilde{\alpha }.\tilde{x}-\widetilde{\alpha }.\tilde{y}\right\|=|\widetilde{\alpha }|\left\|\tilde{x}-\tilde{y}\right\|=\left|\widetilde{\alpha }\right|d(\tilde{x},\tilde{y}).\]

\begin{definition}
Let $(\check{X},\left\|.\right\|)$ be a soft normed linear space and $(Y,A)$ be a non-null member of ${\mathcal S}(\check{X}).$ Then the mapping ${\left\|.\right\|}_Y:SE(Y,A)\to {{\mathcal R}(A)}^*$ given by ${\left\|\tilde{x}\right\|}_Y=\left\|\tilde{x}\right\|$ for all $\tilde{x}\widetilde{\in }(Y,A)$ is a soft norm on $(Y,A)$. This norm ${\left\|.\right\|}_Y$ is known as the relative norm induced on $(Y,A)$ by $\left\|.\right\|$. The soft normed linear space $(Y,{\left\|.\right\|}_Y,A)$ is called a normed subspace or simply a subspace of the soft normed linear space $(\check{X},\left\|.\right\|,A)$.
\end{definition}

\subsection{Sequences and their convergence in soft normed linear spaces}

\begin{definition}
Let $(\check{X},\left\|.\right\|,A)$ be a soft normed linear space and $\tilde{r}\tilde{>}\overline{0}$ be a soft real number. We define the followings;

\begin{equation} \label{eq403} 
B\left(\tilde{x},\tilde{r}\right)=\{\tilde{y}\widetilde{\in }\check{X}:\left\|\tilde{x}-\tilde{y}\right\|\tilde{<}\tilde{r}\}\subset SE(\check{X}) 
\end{equation} 

\begin{equation} \label{eq404} 
\overline{B}\left(\tilde{x},\tilde{r}\right)=\{\tilde{y}\widetilde{\in }\check{X}:\left\|\tilde{x}-\tilde{y}\right\|\widetilde{\le }\tilde{r}\}\subset SE(\check{X}) 
\end{equation} 

\begin{equation} \label{eq405} 
S\left(\tilde{x},\tilde{r}\right)=\{\tilde{y}\widetilde{\in }\check{X}:\left\|\tilde{x}-\tilde{y}\right\|=\tilde{r}\}\subset SE(\check{X}) 
\end{equation} 

$B\left(\tilde{x},\tilde{r}\right),\overline{B}\left(\tilde{x},\tilde{r}\right)$ and $S\left(\tilde{x},\tilde{r}\right)$ are respectively called an open ball, a closed ball and a sphere with centre at $\tilde{x}$ and radius $\tilde{r}$. $SS(B\left(\tilde{x},\tilde{r}\right)),SS(\overline{B}\left(\tilde{x},\tilde{r}\right))$ and $SS(S\left(\tilde{x},\tilde{r}\right))$ are respectively called a soft open ball, a soft closed ball and a soft sphere with centre at $\tilde{x}$ and radius $\tilde{r}$.
\end{definition}

\begin{definition}
A sequence of soft elements $\{{\tilde{x}}_n\}$ in a soft normed linear space $(\check{X},\left\|.\right\|,A)$ is said to be convergent and converges to a soft element $\tilde{x}$ if $\left\|{\tilde{x}}_n-\tilde{x}\right\|\to \overline{0}$ as $n\to \infty $.This means for every $\widetilde{\varepsilon }\tilde{>}\overline{0}$, chosen arbitrarily, $\exists $ a natural number $N=N(\widetilde{\varepsilon })$, such that $\overline{0}\widetilde{\le }\left\|{\tilde{x}}_n-\tilde{x}\right\|\tilde{<}\widetilde{\varepsilon }$ , whenever $n>N$.

\noindent i.e., $n>N\Longrightarrow {\tilde{x}}_n\in B(\tilde{x},\ \widetilde{\varepsilon })$. We denote this by ${\tilde{x}}_n\to \tilde{x}$ as $n\to \infty $ or by ${\mathop{\lim }_{n\to \infty } {\tilde{x}}_n\ }=\tilde{x}$. $\tilde{x}$ is said to be the limit of the sequence ${\tilde{x}}_n$ as $n\to \infty $.
\end{definition} 

\begin{example}
Let us consider the set ${\mathcal R}$ of all real numbers endowed with the usual norm $\left\|.\right\|$. Let $\left(\check{{\mathcal R}},\left\|.\right\|\right)$ or $(\check{{\mathcal R}},\left\|.\right\|,A)$ be the soft norm generated by the crisp norm $\left\|.\right\|,\ $ where $A$ is the non-empty set of parameters. Let $(Y,A)\widetilde{\subset }\check{{\mathcal R}}$ such that $Y\left(\lambda \right)=(0,1]$ in the real line, $\forall \lambda \in A$. Let us choose a sequence $\{{\tilde{x}}_n\}$ of soft elements of $(Y,A)$ where ${\tilde{x}}_n\left(\lambda \right)=\frac{1}{n},\ \forall n\in {\mathcal N},\forall \lambda \in A.$ Then there is no $\tilde{x}\widetilde{\in }(Y,A)$  such that ${\tilde{x}}_n\to \tilde{x}$ in $(Y,{\left\|.\right\|}_Y,A)$. However the sequence $\{{\tilde{y}}_n\}$ of soft elements of $(Y,A)$ where ${\tilde{y}}_n\left(\lambda \right)=\frac{1}{2},\ \forall n\in {\mathcal N},\forall \lambda \in A$ is convergent in $(Y,\left\|.\right\|,A)$ and converges to  $\overline{\frac{1}{2}}$.
\end{example} 

\begin{theorem}
 Limit of a sequence in a soft normed linear space, if exists is unique.
\end{theorem}

\textbf{Proof.}
 If possible let there exists a sequence $\{{\tilde{x}}_n\}$ of soft elements in a soft normed linear space $\left(\check{X},d\right)$ such that ${\mathop{\lim }_{n\to \infty } {\tilde{x}}_n\ }=\tilde{x}$, ${\mathop{\lim }_{n\to \infty } {\tilde{x}}_n\ }={\tilde{x}}^/$, where $\tilde{x}\ne {\tilde{x}}^/$. Then there is at least one $\lambda \in A$ such that 

\noindent $\left\|\tilde{x}-{\tilde{x}}^/\right\|\left(\lambda \right)\ne 0$. We consider a positive real number ${\varepsilon }_{\lambda }$ satisfying $0<{\varepsilon }_{\lambda }<\frac{1}{2}\left\|\tilde{x}-{\tilde{x}}^/\right\|\left(\lambda \right)$. Let $\widetilde{\varepsilon }\tilde{>}\overline{0}$ with $\widetilde{\varepsilon }\left(\lambda \right)={\varepsilon }_{\lambda }$.

\noindent Since ${\tilde{x}}_n\to \tilde{x}$, ${\tilde{x}}_n\to {\tilde{x}}^/$. Corresponding to $\widetilde{\varepsilon }\tilde{>}\overline{0}$, $\exists $ natural numbers $N_1=N_1(\widetilde{\varepsilon })$, $N_2=N_2(\widetilde{\varepsilon })$ such that $n>N_1\Rightarrow {\tilde{x}}_n\in B(\tilde{x},\ \widetilde{\varepsilon })\Rightarrow \left\|\tilde{x}-{\tilde{x}}^/\right\|\tilde{<}\widetilde{\varepsilon }$

\noindent $\Rightarrow \left\|\tilde{x}-{\tilde{x}}^/\right\|\left(\lambda \right)<{\varepsilon }_{\lambda }$, in particular.

\noindent Also, $\ n>N_2\Rightarrow {\tilde{x}}_n\in B({\tilde{x}}^/,\ \widetilde{\varepsilon })\Rightarrow \left\|\tilde{x}-{\tilde{x}}^/\right\|\tilde{<}\widetilde{\varepsilon }\Rightarrow \left\|\tilde{x}-{\tilde{x}}^/\right\|\left(\lambda \right)<{\varepsilon }_{\lambda }$, in particular.

\noindent Hence for all $n>N=\max \{N_1,N_2\}$,

\noindent $\left\|\tilde{x}-{\tilde{x}}^/\right\|\left(\lambda \right)\le \ \left\|{\tilde{x}}_n-\tilde{x}\right\|\left(\lambda \right)+\ \left\|{\tilde{x}}_n-{\tilde{x}}^/\right\|\left(\lambda \right)<2{\varepsilon }_{\lambda }$

\noindent So, ${\varepsilon }_{\lambda }>\frac{1}{2}\left\|\tilde{x}-{\tilde{x}}^/\right\|\left(\lambda \right)$. This contradicts our hypothesis.

\noindent Hence the result follows.

\begin{definition}
A sequence $\{{\tilde{x}}_n\}$ of soft elements in $\left(\check{X},d\right)$ is said to be bounded if the set $\{\left\|{\tilde{x}}_n-{\tilde{x}}_m\right\|;m,n\in N\}$ of soft real numbers is bounded, i.e.,$\exists \widetilde{M}\tilde{>}\overline{0}$ such that $\left\|{\tilde{x}}_n-{\tilde{x}}_m\right\|\widetilde{\le }\widetilde{M}$, $\forall m,n\in N$.
\end{definition}

\begin{definition}
A sequence $\{{\tilde{x}}_n\}$ of soft elements ina soft normed linear space 

\noindent $(\check{X},\left\|.\right\|,A)$ is said to be a Cauchy sequence in $\check{X}$ if corresponding to every $\widetilde{\varepsilon }\tilde{>}\overline{0},\ \exists \ m\in N\ $ such that $\left\|{\tilde{x}}_i-{\tilde{x}}_j\right\|\widetilde{\le }\widetilde{\varepsilon }$, $\forall \ i,j\ge m$  i.e.,  $\left\|{\tilde{x}}_i-{\tilde{x}}_j\right\|\to \overline{0}$ as $i,j\to \infty $.
\end{definition}

\begin{theorem}
Every convergent sequence in a soft normed linear space is Cauchy and every Cauchy sequence is bounded.
\end{theorem}

\textbf{Proof.}
Let $\left\{{\tilde{x}}_n\right\}$ be a convergent sequence of soft elements with limit $\tilde{x}$(say) in $\left(\check{X},\left\|.\right\|\right)$. Then corresponding to each $\widetilde{\varepsilon }\tilde{>}\overline{0},\ \exists \ m\in N\ $ such that ${\tilde{x}}_n\in B\left(\tilde{x},\ \frac{\widetilde{\varepsilon }}{2}\right)$ i.e., $\left\|\tilde{x}-{\tilde{x}}_n\right\|\widetilde{\le }\frac{\widetilde{\varepsilon }}{2}$, $\forall n\ge m$.

\noindent Then for $i,j\ge m$, $\left\|{\tilde{x}}_i-{\tilde{x}}_j\right\|\widetilde{\le }\left\|{\tilde{x}}_i-\tilde{x}\right\|+\left\|\tilde{x}-{\tilde{x}}_j\right\|\tilde{<}\frac{\widetilde{\varepsilon }}{2}+\frac{\widetilde{\varepsilon }}{2}=\widetilde{\varepsilon }$. Hence $\left\{{\tilde{x}}_n\right\}$ is a Cauchy sequence.

\noindent Next let $\left\{{\tilde{x}}_n\right\}$ be a Cauchy sequence of soft elements in $\left(\check{X},\left\|.\right\|\right)$. Then $\exists \ m\in N$ such that  $\left\|{\tilde{x}}_i-{\tilde{x}}_j\right\|\tilde{<}\overline{1}$, $\forall i,j\ge m$. Take $\tilde{M}$ with 

\noindent $\tilde{M}\left(\lambda \right)={\mathop{\max }_{1\le i,j\le m} \left\{\left\|{\tilde{x}}_i-{\tilde{x}}_j\right\|\left(\lambda \right)\right\}\ },\ \forall \lambda \in A$. Then for  $1\le i\le m$ and $j\ge m$, $\left\|{\tilde{x}}_i-{\tilde{x}}_j\right\|\widetilde{\le }\left\|{\tilde{x}}_i-{\tilde{x}}_m\right\|+\left\|{\tilde{x}}_m-{\tilde{x}}_j\right\|\tilde{<}\tilde{M}+\overline{1}$.

\noindent Thus, $\left\|{\tilde{x}}_i-{\tilde{x}}_j\right\|\tilde{<}\tilde{M}+\overline{1}$, $\forall i,j\in N$ and consequently the sequence is bounded.

\begin{definition}
A soft subset $(Y,A)$ with $Y(\lambda )\ne \emptyset $, $\forall \lambda \in A$, in a soft normed linear space $(\check{X},\left\|.\right\|,A)$ is said to be bounded if $\exists $ a soft real number $\tilde{k}$ such that $\left\|\tilde{x}\right\|\widetilde{\le }\tilde{k}$, $\forall \tilde{x}\in (Y,A)$.
\end{definition}

\begin{definition}
Let $(\check{X},\left\|.\right\|,A)$ be a soft normed linear space. Then $\check{X}$ is said to be complete if every Cauchy sequence in $\check{X}$ converges to a soft element of $\check{X}$ i.e., every complete soft normed linear space is called a soft Banach's Space.
\end{definition}

\begin{theorem}
Let $(\check{X},\left\|.\right\|,A)$ be a soft normed linear space. Then

\begin{itemize}

\item [{(i).}] if ${\tilde{x}}_n\to \tilde{x}$ and ${\tilde{y}}_n\to \tilde{y}$ then ${\tilde{x}}_n+{\tilde{y}}_n\to \tilde{x}+\tilde{y}$.

\item [{(ii).}] if ${\tilde{x}}_n\to \tilde{x}$ and ${\widetilde{\lambda }}_n\to \widetilde{\lambda }$ then ${\widetilde{\lambda }}_n.{\tilde{x}}_n\to \widetilde{\lambda }.\tilde{x}$, where $\{{\widetilde{\lambda }}_n\}$ is a sequence of soft scalars.

\item [{(iii).}] if $\{{\tilde{x}}_n\}$ and $\{{\tilde{y}}_n\}$ are Cauchy sequences in $\check{X}$ and ${\{\widetilde{\lambda }}_n\}$ is a Cauchy sequence of soft scalars, then $\{{\tilde{x}}_n+{\tilde{y}}_n\}$ and ${\{\widetilde{\lambda }}_n.{\tilde{x}}_n\}$ are also Cauchy sequences in $\check{X}$.
\end{itemize}
\end{theorem}

\textbf{Proof.}
 (i) Since ${\tilde{x}}_n\to \tilde{x}$ and ${\tilde{y}}_n\to \tilde{y}$, for $\widetilde{\varepsilon }\tilde{>}\overline{0},\ \exists $ +ve integers $N_1,N_2$ such that $\left\|{\tilde{x}}_n-\tilde{x}\right\|\tilde{<}\frac{\widetilde{\varepsilon }}{2},\ \forall n\ge N_1$ and $\left\|{\tilde{y}}_n-\tilde{y}\right\|\tilde{<}\frac{\widetilde{\varepsilon }}{2},\ \forall n\ge N_2.$ Let $N=\max \{N_1,N_2\}$, then both the above relations hold for $n\ge N$.

\noindent Then $\left\|{(\tilde{x}}_n+{\tilde{y}}_n)-(\tilde{x}+\tilde{y})\right\|=\left\|{\tilde{x}}_n+{\tilde{y}}_n-\tilde{x}-\tilde{y}\right\|\widetilde{\le }\left\|{\tilde{x}}_n-\tilde{x}\right\|+\left\|{\tilde{y}}_n-\tilde{y}\right\|\tilde{<}\frac{\widetilde{\varepsilon }}{2}+\frac{\widetilde{\varepsilon }}{2}=\widetilde{\varepsilon },\ \forall n\ge N$.
\[\Longrightarrow {\tilde{x}}_n+{\tilde{y}}_n\to \tilde{x}+\tilde{y}.\]

\noindent (ii) Since ${\tilde{x}}_n\to \tilde{x}$ and ${\lambda }_n\to \lambda $ we get, for $\widetilde{\varepsilon }\tilde{>}\overline{0},\ \exists $ +ve integers $N$ such that $\left\|{\tilde{x}}_n-\tilde{x}\right\|\tilde{<}\widetilde{\varepsilon },\ \forall n\ge N$.

\noindent Now, $\left\|{\tilde{x}}_n\right\|=\left\|{\tilde{x}}_n-\tilde{x}+\tilde{x}\right\|\widetilde{\le }\left\|{\tilde{x}}_n-\tilde{x}\right\|+\left\|\tilde{x}\right\|\tilde{<}\widetilde{\varepsilon }+\left\|\tilde{x}\right\|,\ \forall n\ge N$.

\begin{equation} \label{eq412}
\Longrightarrow \left\|{\tilde{x}}_n\right\|\tilde{<}\widetilde{\varepsilon }+\left\|\tilde{x}\right\|,\ \forall n\ge N.
\end{equation}

\noindent Thus the sequence $\{\left\|{\tilde{x}}_n\right\|\}$ is bounded.

\noindent Now, $\left\|{\widetilde{\lambda }}_n.{\tilde{x}}_n-\widetilde{\lambda }.\tilde{x}\right\|=\left\|{\widetilde{\lambda }}_n.{\tilde{x}}_n-\widetilde{\lambda }.{\tilde{x}}_n+\widetilde{\lambda }.{\tilde{x}}_n-\widetilde{\lambda }.\tilde{x}\right\|=\left\|{{\tilde{x}}_n(\widetilde{\lambda }}_n-\widetilde{\lambda })+\widetilde{\lambda}({\tilde{x}}_n-\tilde{x})\right\|$

\noindent $\widetilde{\le }\left\|{{\tilde{x}}_n(\widetilde{\lambda }}_n-\widetilde{\lambda })\right\|+\left\|\widetilde{\lambda }({\tilde{x}}_n-\tilde{x})\right\|$ 

\noindent $=\left|{\widetilde{\lambda }}_n-\widetilde{\lambda }\right|\left\|{\tilde{x}}_n\right\|+\left|\widetilde{\lambda }\right|\left\|{\tilde{x}}_n-\tilde{x}\right| $

\begin{equation} \label{eq413}
\Longrightarrow \left\|{\widetilde{\lambda }}_n.{\tilde{x}}_n-\widetilde{\lambda }.\tilde{x}\right\|\widetilde{\le}\left|{\widetilde{\lambda }}_n-\widetilde{\lambda }\right|\left\|{\tilde{x}}_n\right\|+\left|\widetilde{\lambda }\right|\left\|{\tilde{x}}_n-\tilde{x}\right|
\end{equation}

Since ${\tilde{x}}_n\to \tilde{x}$ and ${\widetilde{\lambda }}_n\to \widetilde{\lambda }$ we get,$\left|{\widetilde{\lambda }}_n-\widetilde{\lambda }\right|\to \overline{0}$ and $\left\|{\tilde{x}}_n-\tilde{x}\right\|\to \overline{0}$ as $n\to \infty $.

\noindent Now using (\ref{eq412}) and (\ref{eq413}) we get, $\left\|{\widetilde{\lambda }}_n.{\tilde{x}}_n-\widetilde{\lambda }.\tilde{x}\right\|\to \overline{0}$ as $n\to \infty .$ 

\noindent Hence ${\widetilde{\lambda }}_n.{\tilde{x}}_n\to \widetilde{\lambda }.\tilde{x}$.

\noindent 

\noindent (iii) Let $\{{\tilde{x}}_n\}$ and $\{{\tilde{y}}_n\}$ be Cauchy sequences in $\check{X}$, then for $\widetilde{\varepsilon }\tilde{>}\overline{0},\ \exists $ +ve integers $N_1,N_2$ such that 

\noindent $\left\|{\tilde{x}}_n-{\tilde{x}}_m\right\|\tilde{<}\frac{\widetilde{\varepsilon }}{2},\ \forall \ m,n\ge N_1$ and $\left\|{\tilde{y}}_n-{\tilde{y}}_n\right\|\tilde{<}\frac{\widetilde{\varepsilon }}{2},\ \forall \ m,n\ge N_2.$ 

\noindent Let $N=\max \{N_1,N_2\}$, then both the above relations hold for $m,n\ge N$. 

\noindent Now, $\left\|{(\tilde{x}}_n+{\tilde{y}}_n)-{(\tilde{x}}_m+{\tilde{y}}_m)\right\|=\left\|{(\tilde{x}}_n-{\tilde{x}}_m)+({\tilde{y}}_n-{\tilde{y}}_m)\right\|$

\noindent $\widetilde{\le }\left\|{\tilde{x}}_n-{\tilde{x}}_m\right\|+\left\|{\tilde{y}}_n-{\tilde{y}}_m\right\|\tilde{<}\frac{\widetilde{\varepsilon }}{2}+\frac{\widetilde{\varepsilon }}{2}=\widetilde{\varepsilon },\ $

\noindent $\forall \ m,n\ge N$. $\Longrightarrow {\{\tilde{x}}_n+{\tilde{y}}_n\}$ is a Cauchy sequences in $\check{X}$.

\noindent Since ${\{\tilde{x}}_n\}$ is a Cauchy sequences in $\check{X}$, for $\widetilde{\varepsilon }\tilde{>}\overline{0},\ \exists $ +ve integers $N$ such that $\left\|{\tilde{x}}_n-{\tilde{x}}_m\right\|\tilde{<}\widetilde{\varepsilon },\ \forall \ m,n\ge N.$

\noindent Taking in particular $n=m+1$, $\left\|{\tilde{x}}_{m+1}\right\|\tilde{<}\widetilde{\varepsilon },\ \forall \ m,n\ge N,$ so $\{\left\|{\tilde{x}}_n\right\|\}$ is bounded. Now $\{{\widetilde{\lambda }}_n\}$ is bounded too.

\noindent Then, $\left\|{\widetilde{\lambda }}_n.{\tilde{x}}_n-{\lambda }_m.{\tilde{x}}_m\right\|=\left\|{\widetilde{\lambda }}_n.{\tilde{x}}_n-{\widetilde{\lambda }}_n.{\tilde{x}}_m+{\widetilde{\lambda }}_n.{\tilde{x}}_m-{\widetilde{\lambda }}_m.{\tilde{x}}_m\right\|$

\noindent $=\left\|{\widetilde{\lambda }}_n{(\tilde{x}}_n-{\tilde{x}}_m)+{\tilde{x}}_m({\widetilde{\lambda }}_n-{\widetilde{\lambda }}_m)\right\|\widetilde{\le }\left|{\widetilde{\lambda }}_n\right|\left\|{(\tilde{x}}_n-{\tilde{x}}_m)\right\|+\left\|{\tilde{x}}_m\right\|\left|({\widetilde{\lambda }}_n-{\widetilde{\lambda }}_m)\right|\to \overline{0}$ as $n\to \infty $.

\noindent $\Longrightarrow \{{\widetilde{\lambda }}_n.{\tilde{x}}_n\}$ are also Cauchy sequences in $\check{X}$.

\begin{theorem}
If $(M,A)$ is a soft subspace in a soft normed linear space $(\check{X},\left\|.\right\|,A)$, then the closure of $(M,A)$, $\overline{(M,A)}$ is also a soft subspace.
\end{theorem}

\textbf{Proof.}
Let $\tilde{x},\tilde{y}\widetilde{\in }\overline{(M,A)}$, we must show that any linear combination of $\tilde{x},\tilde{y}$ belongs to $\overline{(M,A)}$. Since $\tilde{x},\tilde{y}\widetilde{\in }\overline{(M,A)}$, corresponding to $\widetilde{\varepsilon }\tilde{>}\overline{0}$, there exists soft elements ${\tilde{x}}_1,{\tilde{y}}_1\widetilde{\in }\overline{(M,A)}$ such that $\left\|\tilde{x}-{\tilde{x}}_1\right\|\tilde{<}\widetilde{\varepsilon },\left\|\tilde{y}-{\tilde{y}}_1\right\|\tilde{<}\widetilde{\varepsilon }.$
 
\noindent For soft scalars $\widetilde{\alpha },\widetilde{\beta }\tilde{>}\overline{0}$, $\left\|\left(\widetilde{\alpha }\tilde{x}+\widetilde{\beta }\tilde{y}\right)-(\widetilde{\alpha }{\tilde{x}}_1+\widetilde{\beta }{\tilde{y}}_1)\right\|$

\noindent $\widetilde \le \left|\widetilde{\alpha }\right|\left\|\tilde{x}-{\tilde{x}}_1\right\|+\left|\widetilde{\beta }\right|\left\|\tilde{y}-{\tilde{y}}_1\right\|\tilde{<}\widetilde{\varepsilon }\left(\left|\widetilde{\alpha }\right|+\left|\widetilde{\beta }\right|\right)={\widetilde{\varepsilon }}^/$(say),

\noindent The above inequality shows that $\widetilde{\alpha }{\tilde{x}}_1+\widetilde{\beta }{\tilde{y}}_1$ belongs to the open ball $B(\widetilde{\alpha }\tilde{x}+\widetilde{\beta }\tilde{y},{\widetilde{\varepsilon }}^/)$. As $\widetilde{\alpha }{\tilde{x}}_1+\widetilde{\beta }{\tilde{y}}_1$ and ${\widetilde{\varepsilon }}^/\tilde{>}\overline{0}$ are arbitrary, it follows that $\widetilde{\alpha }\tilde{x}+\widetilde{\beta }\tilde{y}\widetilde{\in }\overline{(M,A)}$. Hence $\overline{(M,A)}$ is a soft subspace of $\check{X}$.

\begin{definition}
 A soft linear space $\check{X}$ is said to be of finite dimensional if there is a finite set of linearly independent soft vectors in $\check{X}$ which also generates $\check{X}$.
\end{definition}

\begin{lemma} \label{lem410}
 Let ${\tilde{x}}_1,{\tilde{x}}_2,\dots ..,{\tilde{x}}_n$ be a linearly independent set of soft vectors in a soft linear space $\check{X}.$ Then there is a soft real number $\tilde{c}\tilde{>}\overline{0}$ such that for every set of soft scalars ${\widetilde{\alpha }}_1,{\widetilde{\alpha }}_2,\dots ..,{\widetilde{\alpha }}_n$ we have
\[\left\|{\widetilde{\alpha }}_1{\tilde{x}}_1+{\widetilde{\alpha }}_2{\tilde{x}}_2+\dots +{\widetilde{\alpha }}_n{\tilde{x}}_n\right\|\widetilde{\ge }\tilde{c}\left(\left|{\widetilde{\alpha }}_1\right|+\left|{\widetilde{\alpha }}_2\right|+\dots +\left|{\widetilde{\alpha }}_n\right|\right).\] 
\end{lemma}

\textbf{Proof.}
 The theorem will be proved if we can prove
\[\left\|{\widetilde{\alpha }}_1{\tilde{x}}_1+{\widetilde{\alpha }}_2{\tilde{x}}_2+\dots +{\widetilde{\alpha }}_n{\tilde{x}}_n\right\|\left(\lambda \right)\ge \left[\tilde{c}\left(\left|{\widetilde{\alpha }}_1\right|+\left|{\widetilde{\alpha }}_2\right|+\dots +\left|{\widetilde{\alpha }}_n\right|\right)\right]\left(\lambda \right),\forall \lambda \in A.\] 
i.e.,${\left\|{\widetilde{\alpha }}_1\left(\lambda \right).{\tilde{x}}_1\left(\lambda \right)+{\widetilde{\alpha }}_2\left(\lambda \right).{\tilde{x}}_2\left(\lambda \right)+\dots +{\widetilde{\alpha }}_n\left(\lambda \right).{\tilde{x}}_n\left(\lambda \right)\right\|}_{\lambda }$

\noindent $\ge \left[\tilde{c}\left(\lambda \right).\left(\left|{\widetilde{\alpha }}_1\left(\lambda \right)\right|+\left|{\widetilde{\alpha }}_2\left(\lambda \right)\right|+\dots +\left|{\widetilde{\alpha }}_n\left(\lambda \right)\right|\right)\right],\forall \lambda \in A.$

\noindent Now, ${\tilde{x}}_1,{\tilde{x}}_2,\dots ..,{\tilde{x}}_n$ being soft vectors in $\check{X},{\tilde{x}}_1\left(\lambda \right),{\tilde{x}}_2\left(\lambda \right),\dots ..,{\tilde{x}}_n\left(\lambda \right)$ are vectors in $X$ and ${\widetilde{\alpha }}_1,{\widetilde{\alpha }}_2,\dots ..,{\widetilde{\alpha }}_n$ being soft scalars ${\widetilde{\alpha }}_1\left(\lambda \right),{\widetilde{\alpha }}_2\left(\lambda \right),\dots ..,{\widetilde{\alpha }}_n\left(\lambda \right)$ are scalars. Then using the property of normed linear space $(X,{\left\|.\right\|}_{\lambda })$ we get a real number $c_{\lambda },$ such that the above relation holds for $\tilde{c}\left(\lambda \right)=c_{\lambda },\forall \lambda \in A.$

\begin{theorem}
 Every Cauchy sequence in ${\mathcal R}(A)$ with finite parametr set $A$ is convergent, i.e., the set of all soft real numbers with its usual modulus soft norm as defined in Example \ref{ex411}, with finite parametr set $A$, is a soft Banach space.
\end{theorem}

\textbf{Proof.}
Let $\{{\tilde{x}}_n\}$ be any arbitrary Cauchy sequence in ${\mathcal R}(A).$ Then corresponding to every $\widetilde{\varepsilon }\tilde{>}\overline{0},\ \exists \ m\in N\ $ such that $|{\tilde{x}}_i-{\tilde{x}}_j|\widetilde{\le }\widetilde{\varepsilon }$, $\forall \ i,j\ge m$ i.e.,$|{\tilde{x}}_i-{\tilde{x}}_j|\left(\lambda \right)\le \widetilde{\varepsilon }\left(\lambda \right)$, $\forall \ i,j\ge m$,$\ \forall \lambda \in A$ i.e.,$|{\tilde{x}}_i\left(\lambda \right)-{\tilde{x}}_j\left(\lambda \right)|\le \widetilde{\varepsilon }\left(\lambda \right)$, $\forall \ i,j\ge m$,$\ \forall \lambda \in A.$ Then $\{{\tilde{x}}_n\left(\lambda \right)\}$ is a Cauchy sequence of ordinary real numbers ${\mathcal R}$ for each $\lambda \in A.$ By the Completeness of ${\mathcal R}$ and finiteness of $A$, it follows that $\{{\tilde{x}}_n\left(\lambda \right)\}$ is convergent for each $\lambda \in A.$ Let ${\tilde{x}}_n\left(\lambda \right)\to x_{\lambda },$ for each $\lambda \in A.$ Consider the soft element $\tilde{x}$ defined by $\tilde{x}\left(\lambda \right)=x_{\lambda },$ for each $\lambda \in A.$ Then $\tilde{x}$ is a soft real number and it follows that the sequence $\{{\tilde{x}}_n\}$ of soft real numbers is convergent and it converges to the soft real number $\tilde{x}.$ Hence ${\mathcal R}(A)$ is a soft Banach space.

\begin{theorem}
 Every finite dimensional soft normed linear space over a finite parameter set $A$ is complete.
\end{theorem}

\textbf{Proof.}
Let $\check{X}$ be a finite dimensional soft normed linear space over a finite parameter set $A$. Let $\{{\tilde{y}}_m\}$ be any arbitrary Cauchy sequence in $\check{X}.$ We show that $\{{\tilde{y}}_m\}$ converges to some soft element $\tilde{y}\widetilde{\in }\check{X}.$ Suppose that the dimension of $\check{X}$ is $n,$ and let $\{{\tilde{e}}_1,{\tilde{e}}_2,\dots ..,{\tilde{e}}_n\}$ be a basis for $\check{X}.$ Then each ${\tilde{y}}_m$ has a unique representation ${\tilde{y}}_m={{\widetilde{\alpha }}_1}^{(m)}{\tilde{e}}_1+{{\widetilde{\alpha }}_2}^{(m)}{\tilde{e}}_2+\dots +{{\widetilde{\alpha }}_n}^{(m)}{\tilde{e}}_n.$

\noindent Because $\{{\tilde{y}}_m\}$ is a Cauchy sequence, for $\widetilde{\varepsilon }\tilde{>}\overline{0}$ arbitrary there exist a positive integer $N$ such that $\left\|{\tilde{y}}_m-{\tilde{y}}_r\right\|\tilde{<}\widetilde{\varepsilon }$  for $m,r>N.$

\noindent From Lemma \ref{lem410}, it follows that there exists $\tilde{c}\tilde{>}\overline{0}$ such that

\noindent $\widetilde{\varepsilon }\tilde{>}\left\|{\tilde{y}}_m-{\tilde{y}}_r\right\|=\left\|\sum^n_{j=1}{({{\widetilde{\alpha }}_j}^{\left(m\right)}-{{\widetilde{\alpha }}_j}^{(r)}){\tilde{e}}_j}\right\|\widetilde{\ge }\tilde{c}\sum^n_{j=1}{\left|{{\widetilde{\alpha }}_j}^{\left(m\right)}-{{\widetilde{\alpha }}_j}^{\left(r\right)}\right|,}$ for $m,r>N.$

\noindent Consequently, $\left|{{\widetilde{\alpha }}_j}^{\left(m\right)}-{{\widetilde{\alpha }}_j}^{\left(r\right)}\right|\widetilde{\le }\sum^n_{j=1}{\left|{{\widetilde{\alpha }}_j}^{\left(m\right)}-{{\widetilde{\alpha }}_j}^{\left(r\right)}\right|}\tilde{<}{\widetilde{\varepsilon }}/{\tilde{c}}$

\noindent shows that each of the $n$ sequences $\left\{{{\widetilde{\alpha }}_j}^{\left(m\right)}\right\}=\left\{{{\widetilde{\alpha }}_j}^{\left(1\right)},{{\widetilde{\alpha }}_j}^{\left(2\right)},{{\widetilde{\alpha }}_j}^{\left(3\right)},\dots ..\right\},\ j=1,2,..,n$ is Cauchy in ${\mathcal R}(A)$ and $A$ is finite,  converges to ${\widetilde{\alpha }}_j,$ (say), $j=1,2,\dots ,n.$

\noindent We now define the soft element $\tilde{y}={\widetilde{\alpha }}_1{\tilde{e}}_1+{\widetilde{\alpha }}_2{\tilde{e}}_2+\dots +{\widetilde{\alpha }}_n{\tilde{e}}_n$ which is clearly a soft element of $\check{X}.$ Moreover, since ${{\widetilde{\alpha }}_j}^{\left(m\right)}\to {\widetilde{\alpha }}_j$ as $m\to \infty $ and $j=1,2,\dots ,n;$ we have 

\noindent $\left\|{\tilde{y}}_m-\tilde{y}\right\|=\left\|\sum^n_{j=1}{({{\widetilde{\alpha }}_j}^{\left(m\right)}-{\widetilde{\alpha }}_j){\tilde{e}}_j}\right\|\widetilde{\le }\sum^n_{j=1}{\left|{{\widetilde{\alpha }}_j}^{\left(m\right)}-{\widetilde{\alpha }}_j\right|}$$\left\|{\tilde{e}}_j\right\|\to \overline{0}$ as $m\to \infty .$

\noindent i.e., ${\tilde{y}}_m\to \tilde{y}$ as $m\to \infty .$

\subsection{Equivalent soft norms}

\begin{definition}
 Let $\check{X}$ be a soft linear (vector) space. A soft norm ${\left\|.\right\|}_1$ on $\check{X}$ is said to be equivalent to a soft norm ${\left\|.\right\|}_2$ on $\check{X}$ if there are positive soft real numbers $\tilde{a}$ and $\tilde{b}$ such that for all $\tilde{x}\widetilde{\in }\check{X}$ we have 

\begin{equation} \label{eq431}
\tilde{a}{\left\|\tilde{x}\right\|}_2\widetilde{\le }{\left\|\tilde{x}\right\|}_1\widetilde{\le }\tilde{b}{\left\|\tilde{x}\right\|}_2
\end{equation}
\end{definition} 

\begin{theorem}
 On a finite dimensional soft linear space $\check{X}$, any soft norm ${\left\|.\right\|}_1$ is equivalent to any other norm ${\left\|.\right\|}_2.$
\end{theorem}

\textbf{Proof.}
 Let $n$ be the dimension of $\check{X}$ and $\{{\tilde{e}}_1,{\tilde{e}}_2,\dots ..,{\tilde{e}}_n\}$ be a basis for $\check{X}.$ If $\tilde{x}\widetilde{\in }\check{X},$ then $\tilde{x}$ has the representation $\tilde{x}={{\widetilde{\alpha }}_1\tilde{e}}_1+{{\widetilde{\alpha }}_2\tilde{e}}_2+\dots .+{{\widetilde{\alpha }}_n\tilde{e}}_n$.

\noindent By Lemma \ref{lem410}, there is a soft real number $\tilde{c}\tilde{>}\overline{0}$ such that, 

\noindent ${\left\|\tilde{x}\right\|}_1\widetilde{\ge }\tilde{c}\left(\left|{\widetilde{\alpha }}_1\right|+\left|{\widetilde{\alpha }}_2\right|+\dots +\left|{\widetilde{\alpha }}_n\right|\right)$.

\noindent If $\tilde{R}\left(\lambda \right)={\mathop{\max }_{j} \{{\left\|{\tilde{e}}_j\right\|}_2\ }\left(\lambda \right)\},\forall \lambda \in A.$ Then soft norm axioms give,
\[{\left\|\tilde{x}\right\|}_2\widetilde{\le }\sum^n_{j=1}{\left|{\widetilde{\alpha }}_j\right|}{\left\|{\tilde{e}}_j\right\|}_2\widetilde{\le }\tilde{R}.\sum^n_{j=1}{\left|{\widetilde{\alpha }}_j\right|}\widetilde{\le }\left(\frac{\tilde{R}}{\tilde{c}}\right).{\left\|\tilde{x}\right\|}_1\] 
or, $\left(\frac{\tilde{c}}{\tilde{R}}\right).{\left\|\tilde{x}\right\|}_2\widetilde{\le }{\left\|\tilde{x}\right\|}_1$

\noindent The other side inequality in (\ref{eq431}) is obtained by interchanging the roles of ${\left\|.\right\|}_1$ and ${\left\|.\right\|}_2$ in the above argument.

\begin{lemma}{(Riesz's Lemma)}
 Let $\tilde{L}$ be a proper soft closed subspace of a soft normed linear space $\check{X}$ satisfying (N5). Then for $\widetilde{\varepsilon }\tilde{>}\overline{0},\ $ there exists $\tilde{y}\widetilde{\in }\check{X}-\tilde{L}$ with $\left\|\tilde{y}\right\|=\overline{1}$ such that for all $\tilde{x}\widetilde{\in }\tilde{L}$, the inequality $\left\|\tilde{x}-\tilde{y}\right\|\tilde{>}\overline{1}-\widetilde{\varepsilon }$ is satisfied.
\end{lemma} 

\textbf{Proof.}
Let $\widetilde{\varepsilon }\tilde{>}\overline{0},$ then $\widetilde{\varepsilon }\left(\lambda \right)={\varepsilon }_{\lambda }>0,\forall \lambda \in A.$

\noindent Also, since $\check{X}$ satisfies (N5),$\tilde{L}\left(\lambda \right)=L_{\lambda }$ is a proper closed subspace of the normed linear space $\check{X}\left(\lambda \right)=X,$ for each $\lambda \in A.$ Thus by Riesz's Lemma for normed linear space $(X,{\left\|.\right\|}_{\lambda })$, there exists $\tilde{y}\left(\lambda \right)\in X-L_{\lambda }$ with ${\left\|\tilde{y}\left(\lambda \right)\right\|}_{\lambda }=1$ such that for all $\tilde{x}\left(\lambda \right)\in L_{\lambda }$ the inequality ${\left\|\tilde{x}\left(\lambda \right)-\tilde{y}\left(\lambda \right)\right\|}_{\lambda }>1-{\varepsilon }_{\lambda }$ is satisfied.

\noindent Then for $\widetilde{\varepsilon }\tilde{>}\overline{0},$ there exists  $\tilde{y}\widetilde{\in }\check{X}-\tilde{L}$ with $\left\|\tilde{y}\right\|=\overline{1}$ such that for all $\tilde{x}\widetilde{\in }\tilde{L}$, the inequality $\left\|\tilde{x}-\tilde{y}\right\|\tilde{>}\overline{1}-\widetilde{\varepsilon }$ is satisfied.

\subsection{Convex sets in soft normed linear spaces}

\begin{definition}
 Let $\check{X}$ be a soft normed linear space and ${\tilde{x}}_1,\ {\tilde{x}}_2\widetilde{\in }\check{X}.$ The set of all soft elements of the form $\tilde{y}=\tilde{t}{\tilde{x}}_1+\left(\overline{1}-\tilde{t}\right){\tilde{x}}_2,$ where $\tilde{t}$ assumes all soft real numbers such that $\tilde{t}\left(\lambda \right)\in \left[0,\ 1\right],\ \forall \lambda \in A,$ is called the segment joining the soft elements ${\tilde{x}}_1$ and ${\tilde{x}}_2$. A soft set $\tilde{K}\widetilde{\subset }\check{X}$ is called convex if all segments joining any two soft elements of $\tilde{K}$ are contained in $\tilde{K}$. Clearly, every soft subspace $\tilde{M}$ of $\check{X}$ is a convex soft set.
\end{definition}

\begin{theorem} 
A soft sphere in a soft normed linear space is a convex soft set.
\end{theorem} 

\textbf{Proof.}
 We prove the theorem for a soft closed sphere. The proof for the soft open sphere is analogous. Let ${\tilde{x}}_1,\ {\tilde{x}}_2\widetilde{\in }SS\left(S(\tilde{a},\tilde{r})\right),$ so that $\left\|{\tilde{x}}_1-\tilde{a}\right\|\widetilde{\le }\tilde{r}$ and $\left\|{\tilde{x}}_2-\tilde{a}\right\|\widetilde{\le }\tilde{r}$. Let $\tilde{y}=\tilde{t}{\tilde{x}}_1+\left(\overline{1}-\tilde{t}\right){\tilde{x}}_2,\ \tilde{t}\left(\lambda \right)\in \left[0,\ 1\right],\ \forall \lambda \in A.$ Then we have, $\left\|\tilde{y}-\tilde{a}\right\|=\left\|\tilde{t}{\tilde{x}}_1+\left(\overline{1}-\tilde{t}\right){\tilde{x}}_2-\tilde{a}\right\|=\left\|\tilde{t}{\tilde{x}}_1+\left(\overline{1}-\tilde{t}\right){\tilde{x}}_2-\tilde{t}\tilde{a}-\left(\overline{1}-\tilde{t}\right)\tilde{a}\right\|$

\noindent $=\left\|\tilde{t}{\tilde{x}}_1-\tilde{t}\tilde{a}\right\|+\left\|\left(\overline{1}-\tilde{t}\right){\tilde{x}}_2-\left(\overline{1}-\tilde{t}\right)\tilde{a}\right\|=\tilde{t}\left\|{\tilde{x}}_1-\tilde{a}\right\|+\left(\overline{1}-\tilde{t}\right)\left\|{\tilde{x}}_2-\tilde{a}\right\|$

\noindent $\widetilde{\le }\tilde{t}\tilde{r}+\left(\overline{1}-\tilde{t}\right)\tilde{r}=\tilde{r},$ 

\noindent So, $\tilde{y}\widetilde{\in }SS\left(S(\tilde{a},\tilde{r})\right)$. This proves the theorem.

\begin{theorem}
 Let $\check{X}$ be a soft normed linear space and $\tilde{K}$ be a convex soft subset of $\check{X}$. Then the closure of $\tilde{K}$, $\overline{\tilde{K}}$ is convex.
\end{theorem} 

\textbf{Proof.}
Let $\tilde{x},\tilde{y}\widetilde{\in }\overline{\tilde{K}}$ and $\widetilde{\varepsilon }\tilde{>}\overline{0},$ there exists ${\tilde{x}}_1,{\tilde{y}}_1\widetilde{\in }\tilde{K}$ such that 

\noindent $\left\|\tilde{x}-{\tilde{x}}_1\right\|\tilde{<}\widetilde{\varepsilon },\left\|\tilde{y}-{\tilde{y}}_1\right\|\tilde{<}\widetilde{\varepsilon }$.

\noindent Let $\tilde{t}\left(\lambda \right)\in \left[0,\ 1\right],\ \forall \lambda \in A.$ Then $\tilde{t}$ is a non-negative soft real number.

\noindent Then, $\left\|\tilde{t}\tilde{x}+\left(\overline{1}-\tilde{t}\right)\tilde{y}-\left\{\tilde{t}{\tilde{x}}_1+\left(\overline{1}-\tilde{t}\right){\tilde{y}}_1\right\}\right\|\widetilde{\le }\tilde{t}\left\|\tilde{x}-{\tilde{x}}_1\right\|+\left(\overline{1}-\tilde{t}\right)\left\|\tilde{y}-{\tilde{y}}_1\right\|$

\noindent $\widetilde{\le }\tilde{t}\widetilde{\varepsilon }+\left(\overline{1}-\tilde{t}\right)\widetilde{\varepsilon }=\widetilde{\varepsilon };$

\noindent Since $\tilde{K}$ is convex, $\tilde{t}{\tilde{x}}_1+\left(\overline{1}-\tilde{t}\right){\tilde{y}}_1\widetilde{\in }\tilde{K}$ and because $\widetilde{\varepsilon }\tilde{>}\overline{0},$ is arbitrary, 

\noindent $\tilde{t}\tilde{x}+\left(\overline{1}-\tilde{t}\right)\tilde{y}\widetilde{\in }\overline{\tilde{K}}$. $\Rightarrow \overline{\tilde{K}}$ is convex.

\begin{theorem}
 The intersection of an arbitrary number of convex soft sets is a convex soft set.
\end{theorem} 

\textbf{Proof.}
 Let $\tilde{M}=\bigcap_a{{\tilde{M}}_a}$, where each ${\tilde{M}}_a$ is a convex soft set. If $\tilde{x},\tilde{y}\widetilde{\in }\tilde{M}$, then $\tilde{x},\tilde{y}$ belongs to all ${\tilde{M}}_a$ and because each ${\tilde{M}}_a$ is a convex, $\tilde{t}\tilde{x}+\left(\overline{1}-\tilde{t}\right)\tilde{y}\widetilde{\in }{\tilde{M}}_a,\ $ where $\tilde{t}\left(\lambda \right)\in \left[0,\ 1\right],\ \forall \lambda \in A.$ So, $\tilde{t}\tilde{x}+\left(\overline{1}-\tilde{t}\right)\tilde{y}\widetilde{\in }\tilde{M}$ and $\tilde{M}$ is convex.

\end{document}